\documentclass[a4paper]{IEEEtran}
\usepackage{amsmath}
\usepackage{algorithm,algorithmic}
\usepackage[pdftex]{graphicx}
\usepackage{amsmath,amssymb}
\usepackage{bm}
\usepackage{color}

\newtheorem{lemma}{Lemma}
\newtheorem{corollary}{Corollary}
\newtheorem{theorem}{Theorem}

\title{Stochastic Dynamics of Noisy Average Consensus: Analysis and Optimization}

\author{%
  \IEEEauthorblockN{
  		Tadashi Wadayama and Ayano Nakai-Kasai} \\
  \IEEEauthorblockA{\IEEEauthorrefmark{1}%
		Nagoya Institute of Technology,
		Gokiso, Nagoya, Aichi 466-8555, Japan,\\
 		wadayama@nitech.ac.jp, nakai.ayano@nitech.ac.jp}
 \thanks{
Part of this research was presented 
at IEEE International Symposium on Information Theory 2022 (ISIT2022)~\cite{Wadayama2022}.
}

}
\begin{document}

\maketitle

\begin{abstract}
A continuous-time average consensus system
is a linear dynamical system defined over a graph, where
each node has its own state value that evolves 
according to a simultaneous linear differential equation.
A node is allowed to interact with neighboring nodes.
Average consensus is a phenomenon that the all the state 
values converge to the average of the initial state values. 
In this paper, 
we assume that a node can communicate with neighboring nodes
through an additive white Gaussian noise channel. 
We first formulate the noisy average consensus system 
by using a stochastic differential equation (SDE), 
which allows us to use the Euler-Maruyama method, 
a numerical technique for solving SDEs. 
By studying the stochastic behavior of the residual error of the Euler-Maruyama method, 
we arrive at the covariance evolution equation.
The analysis of the residual error 
leads to a compact formula for mean squared error (MSE), 
which shows that the sum of the inverse eigenvalues of the Laplacian matrix 
is the most dominant factor influencing the MSE.
Furthermore, we propose optimization problems aimed at minimizing the MSE 
at a given target time, 
and introduce a deep unfolding-based optimization method 
to solve these problems. 
The quality of the solution is validated by numerical experiments.
\end{abstract}

\section{Introduction}

Continuous-time {\em average consensus system}  
is a linear dynamical system defined 
over a graph \cite{Olfati-Saber07}.
Each node has its own state value, and it evolves 
according to a simultaneous linear differential equation
where a node  is {only allowed }
to interact with neighboring nodes.
The ordinary differential equation (ODE) at the node $i (1 \le i \le n)$ 
governing the evolution of the state value $x_i(t)$ of the node $i$
is given by
\begin{align}
	\frac{dx_i(t)}{dt} = -\sum_{j \in {\cal N}_i} \mu_{ij}(x_i(t) - x_j(t)).
\end{align}
The set ${\cal N}_i$ denote the neighboring nodes of node $i$, 
while the positive scalar $\mu_{ij}$ denotes the edge weight associated with the edge $(i,j)$. 
The same ODE applies to all other nodes as well. 
These dynamics gradually decrease the differences
between the state values of neighboring nodes, 
leading to a phenomenon called average consensus 
that the all the state 
values converge to the average of the initial state values \cite{OlfatiSaber-Murray}. 

 
The average consensus system has been studied in numerous fields 
such as
multi-agent control \cite{Reb},
distributed algorithm \cite{Xiao}, 
formation control \cite{Ren}.
An excellent survey on average consensus systems
can be found in \cite{Olfati-Saber07}.


In this paper, we will examine average consensus systems within the context of 
communications across noisy channels, such as wireless networks. 
Specifically, we consider the scenario in which nodes engage in local wireless communication, 
such as drones flying in the air or sensors dispersed across a designated area. 
It is assumed that each node can only communicate with neighboring nodes 
via an additive white Gaussian noise (AWGN) channel. 
The objective of the communication is to aggregate 
the information held by all nodes through the application 
of average consensus systems. 
As previously stated, the consensus value is the average of the initial state values.

In this setting, we must account for the impact of Gaussian noise on the differential equations. The differential equation for a {\em noisy average consensus system} takes the form:
\begin{align}\label{white_noise_model}
 \frac{d x_i(t)}{dt} = -\sum_{j \in {\cal N}_i} \mu_{ij}(x_i(t) - x_j(t)) + \alpha W_i(t),
\end{align}
where $W_i(t)$ represents an additive white Gaussian process, 
and $\alpha$ is a positive constant. 
The noise $W_i(t)$ can be considered as the sum of the noises occurring on the edges
adjacent to the node $i$.
In a noiseless average consensus system, 
it is well-established that the second smallest eigenvalue of 
the Laplacian matrix of the graph determines 
the convergence speed to the average~\cite{Xiao}. 
The convergence behavior of a noisy system may be quite different from 
that of the noiseless system due to the presence of edge noise. 
However, the stochastic dynamics of such a system has not yet been studied. 
Studies on discrete-time consensus protocols subject to additive noise can be found 
in \cite{Jadbabaie}\cite{Rajagopal}, but to the best of our knowledge, 
there are no prior studies on continuous-time noisy consensus systems.

The main goal of this paper is to study
the stochastic dynamics 
of continuous-time noisy average consensus system.
The theoretical understanding of the stochastic behavior of such systems will be valuable 
for various areas such as multi-agent control and the design of consensus-based distributed algorithms 
for noisy environments.


The primary contributions of this paper are as follows.
We first formulate the noisy average consensus systems 
using stochastic differential equations (SDE)~\cite{Kloeden}\cite{Oksendal}.
This SDE formulation facilitates mathematically rigorous treatment of noisy average consensus.
We use the Euler-Maruyama method~\cite{Kloeden} for solving the SDE, 
which is a numerical method for solving SDEs.
We derive a closed-form mean squared error (MSE) formula 
by analyzing the stochastic behavior of the residual errors in the Euler-Maruyama method. 
We show that the MSE is dominated by the sum of the inverse eigenvalues of the Laplacian matrix. 
However, minimizing the MSE at a specific target time is a non-trivial task because
the objective function involves the sum of the inverse eigenvalues.
To solve this optimization problem, we will propose a deep unfolding-based optimization method.

The outline of the paper is as follows.
In Section 2, we introduce the mathematical notation used throughout the paper, and then provide the definition and fundamental properties of average consensus systems.
In Section 3, we define a noisy average consensus system as a SDE.
In Section 4, we present an analysis of the stochastic behavior of the consensus error and derive a concise MSE formula.
In Section 5, we propose a deep unfolding-based optimization method for minimizing the MSE at a specified target time.
Finally, in Section 6, we conclude the discussion.
\section{Preliminaries} 
\label{sec:preliminaries}

\subsection{Notation}
The following notation will be used throughout this paper.
The symbols $\mathbb{R}$ and $\mathbb{R}_+$ represent 
the set of real numbers and the set of positive real
numbers, respectively.
The one dimensional Gaussian distribution with mean $\mu$ and variance $\sigma^2$ is denoted by ${\cal N}(\mu, \sigma^2)$. The multivariate Gaussian distribution 
with mean vector $\bm \mu$ and covariance matrix
$\bm \Sigma$ is represented by ${\cal N}(\bm \mu, \bm \Sigma)$. 
The expectation operator is denoted by ${\sf E}[\cdot]$.
The notation $\mbox{diag}(\bm x)$ is
the diagonal matrix whose diagonal elements 
are given by $\bm x \in \mathbb{R}^n$.
The matrix exponential $\exp(\bm X) (\bm X \in \mathbb{R}^{n \times n})$ is defined by
\begin{align}
\exp(\bm X) \equiv \sum_{k=0}^\infty \frac{1}{k!}\bm X^k.	
\end{align}
The Frobenius norm of $\bm X\in \mathbb{R}^{n \times n}$
is denoted by $\|\bm X\|_F$.
The notation $[n]$ denotes the set of consecutive integers
from $1$ to $n$. 
  
\subsection{Average Consensus}

Let $G \equiv (V, E)$ be a connected undirected graph where
$V = [n]$. 
Suppose that a node $i \in V$ can be 
regarded as an {\em agent} communicating over the graph $G$.
Namely, a node $i$ and a node $j$ can communicate 
with each other if $(i,j) \in E$. We will not distinguish $(i,j)$ and 
$(j,i)$ because the graph $G$ is undirected.

Each node $i$ has a state value $x_i(t) \in \mathbb{R}$
where $t \in \mathbb{R}$ 
represents continuous-time variable.
The neighborhood of a node $i \in V$ is represented 
by
\begin{align}
{\cal N}_i \equiv \{j \in V: (j,i) \in E, i \ne j\}.	
\end{align}
Note that the node $i$ is excluded from ${\cal N}_i$.
For any time $t$, a node $i \in V$ can access 
the self-state $x_i(t)$ and 
the state values of its neighborhood, i.e.,
$x_j(t), j \in {\cal N}_i$ but 
cannot access to the other state values.

In this section, we briefly review the basic properties of the average consensus processes \cite{Olfati-Saber07}.
We now assume that the set of state values 
$\bm x(t) \equiv (x_1(t), x_2(t), \ldots, x_n(t))^{T}$ are evolved 
according to the simultaneous differential equations
\begin{align} \label{consensus}
	\frac{dx_i(t)}{dt} = -\sum_{j \in {\cal N}_i} \mu_{ij}(x_i(t) - x_j(t)),\quad i \in [n],
\end{align}
where the initial condition is 
\begin{align}
\bm x(0) = \bm c \equiv (c_1, c_2,\ldots, c_n)^{T} \in \mathbb{R}^n.
\end{align}
The edge weight $\mu_{ij}$ follows the symmetric condition
\begin{align}
	\mu_{ij} = \mu_{ji},\quad i\in [n], j\in [n].
\end{align}
Let $\bm \Delta \equiv (\Delta_1,\Delta_2,\ldots,\Delta_n)^T \in \mathbb{R}_{+}^n$ be a {\em degree sequence}
where $\Delta_i$ is defined by
\begin{align} \label{mu_sum}
\Delta_i \equiv \sum_{j \in {\cal N}_i} \mu_{ij},\quad i \in [n].
\end{align}

The continuous-time dynamical system (\ref{consensus}) is called 
an {\em average consensus system} because 
a state value converges to the average of the initial state values at the limit of $t \rightarrow \infty$, i.e, 
\begin{align} \label{average_convergence}
	\lim_{t \rightarrow \infty}\bm x(t) 
	= \frac{1}{n}\left(\sum_{i =1}^n c_i \right) \bm 1 { = \gamma \bm 1},
\end{align}
where the vector $\bm 1$ represents $(1,1,\ldots,1)^{T}$
 and $\gamma$ is defined by 
\begin{align}
\gamma \equiv \frac 1 n \sum_{i =1}^n c_i.	
\end{align}

We define the Laplacian matrix $\bm L \equiv \{L_{ij}\}\in \mathbb{R}^{n \times n}$ of this consensus system as follows:
\begin{align}
	L_{ij} &= \Delta_i,\quad i = j,\ i \in [n], \\
	L_{ij} &= -\mu_{ij},\quad i \ne j \mbox{ and } (i,j) \in E,  \\
	L_{ij} &= 0,\quad i \ne j \mbox{ and } (i,j) \notin E.
\end{align}
From this definition, a Laplacian matrix satisfies
\begin{align}
	\bm L\bm 1 &= \bm 0, \\
	\mbox{diag}(\bm L) &= \bm \Delta, \\
	\bm L &= \bm L^T.
\end{align}

Note that the eigenvalues of the Laplacian matrix 
$\bm{L}$ are { nonnegative real because $\bm L$ is a positive  semi-definite symmetric matrix.}
Let $\lambda_1 = 0 < \lambda_2 \le 
\ldots \le \lambda_n$ be the eigenvalues of $\bm L$
and $\bm \xi_1, \bm \xi_2,\ldots, \bm \xi_n$ be 
the corresponding orthonormal eigenvectors.
The first eigenvector $\bm \xi_1 \equiv (1/\sqrt{n})\bm 1$
is corresponding to the eigenvalue $\lambda_1 = 0$,
which results in $\bm L\bm \xi_1 = 0$.

By using the notion of the Laplacian matrix,
the dynamical system (\ref{consensus}) can be compactly 
rewritten as
\begin{align} \label{consensus2}
	\frac{d\bm x(t)}{dt} = - \bm L \bm x(t),
\end{align}
where the initial condition is 
$\bm x(0) = \bm c$.
The dynamical behaviors of {the average consensus system} (\ref{consensus2})
are thus characterized by the Laplacian matrix $\bm L$.
Since the ODE (\ref{consensus2}) is a linear ODE, it can be easily solved.
The solution of the ODE (\ref{consensus2}) is given by
\begin{align} \label{cont_sol}
	\bm x(t) = \exp(-\bm L t) \bm x(0),\quad t \ge 0.
\end{align}

Let $\bm U \equiv (\bm \xi_1, \bm \xi_2, \ldots, \bm \xi_n) \in \mathbb{R}^{n \times n}$ where $\bm U$ is 
an orthogonal matrix.
The Laplacian matrix $\bm L$ can be diagonalized by using  
$\bm U$, i.e., 
\begin{align}
\bm L = \bm U \mbox{diag}(\lambda_1,\ldots, \lambda_n) \bm U^{T}.
\end{align}

On the basis of the diagonalization, we have the spectral expansion of the matrix exponential:
\begin{align} \nonumber
	\exp(-\bm L t) 
	&= \exp(-\bm U \mbox{diag}(\lambda_1,\ldots, \lambda_n) \bm U^{T} t) \\ \nonumber
	&= \bm U \exp(-\mbox{diag}(\lambda_1,\ldots, \lambda_n)t) \bm U^{T}  \\
	&= \sum_{i=1}^n \exp(-\lambda_i t)\bm \xi_i \bm \xi_i^{T}. 
\end{align}
Substituting this to $\bm x(t) = \exp(-\bm L t) \bm x(0)$, we immediately have
\begin{align}
\bm x(t) = \frac 1 n \bm 1 (\bm 1^{T}) \bm c + \sum_{i=2}^n \exp(-\lambda_i t)\bm \xi_i \bm \xi_i^{T} \bm c.
\end{align}
The second term of the right-hand side converges 
to zero since $\lambda_k > 0$ for $k = 2,3,\ldots, n$.
This explains why average consensus happens, i.e., 
the convergence to the average of the initial state values (\ref{average_convergence}).
The second smallest eigenvalue $\lambda_2$,
called {\em algebraic connectivity} \cite{Chung},  
determines the convergence speed because
$\exp(-\lambda_2 t)\bm \xi_2  \bm \xi_2^{\sf T}$ 
shows the slowest convergence in the second term.

\section{Noisy average consensus system}
\label{sec:noisyconsensus}

\subsection{SDE formulation}

The dynamical model (\ref{white_noise_model}) containing a white Gaussian noise process is mathematically challenging to handle. We will use a common approach of approximating the white Gaussian process by using the standard Wiener process. Instead of model (\ref{white_noise_model}), we will focus on the following stochastic differential equation (SDE) ~\cite{Oksendal}
\begin{equation}\label{sde}
d \bm x(t) = - \bm L \bm x(t) dt + \alpha d \bm b(t)
\end{equation}
to study the noisy average consensus system.
The parameter $\alpha$ is a positive real number, 
and it represents the {\em intensity of the noises}.
The stochastic term $\bm b(t)$ represents the $n$-dimensional standard Wiener process.
The elements of $\bm b(t) = (b_1(t), b_2(t), \ldots, b_n(t))^{T}$ are 
independent one dimensional-standard Wiener processes.
For the Wiener process $b(t)$, we have
$b(0) = 0$, $E[b(t)] = 0$, and it satisfies 
\begin{align}
b(t) - b(s) \sim {\cal N}(0, t-s), \ 0 \le s \le t. 
\end{align}

\subsection{Approaches for studying stochastic dynamics}

%

Our primary objective in the following analysis is to investigate the stochastic dynamics of the noisy average consensus system, focusing on deriving the mean and covariance of the solution $\bm x(t)$ for the SDE (\ref{sde}).

\begin{figure}[htbp]
  \centering
  \includegraphics[width=0.92\hsize]{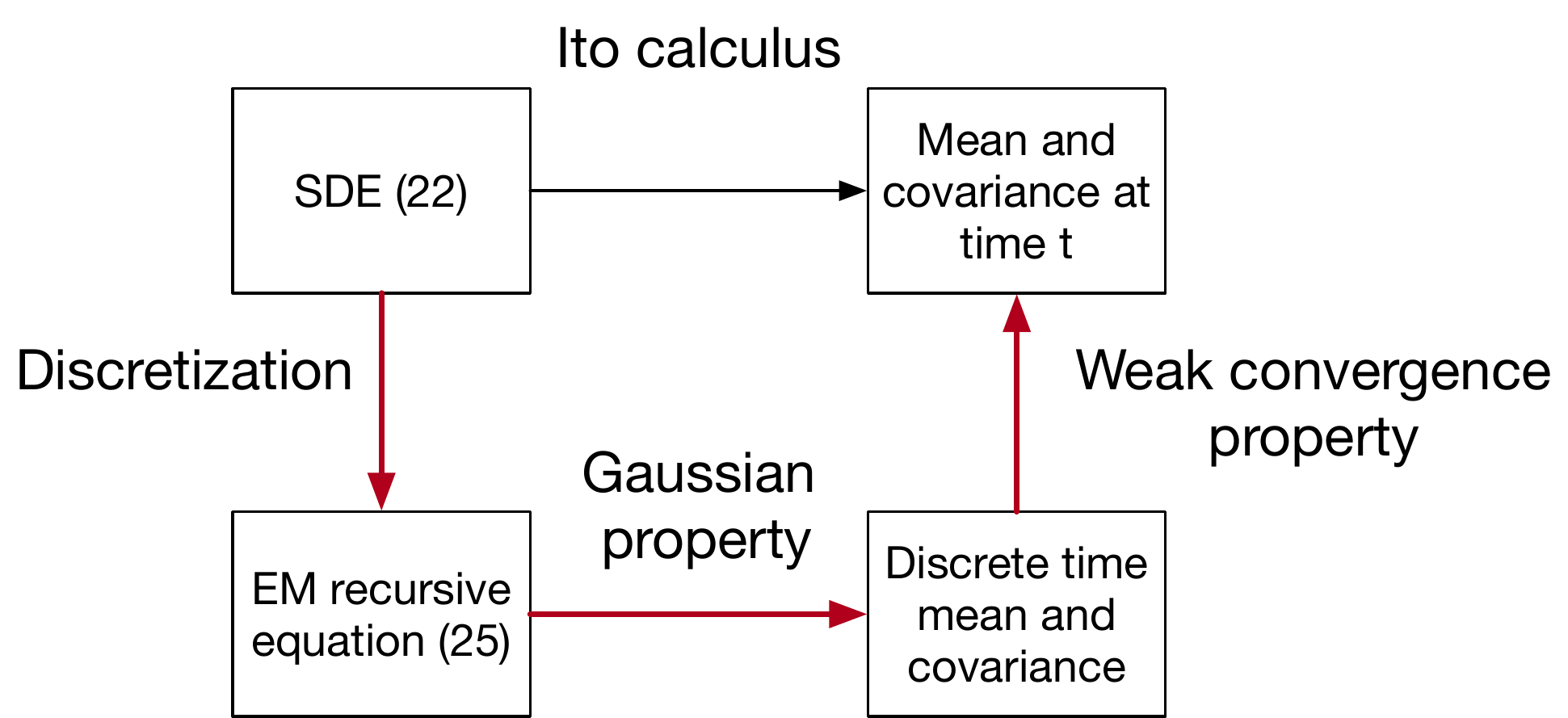}
  \caption{Two approaches for deriving the mean and covariance of $\bm x(t)$. This paper follows the lower path using the EM method.}
  \label{ito}
\end{figure}

There are two approaches to analyze the system. The first approach relies on the established theory of {\em Ito calculus} \cite{Oksendal}, which is used to handle stochastic integrals directly (see Fig.~\ref{ito}). Ito calculus can be applied to derive the first and second moments of the solution of (\ref{sde}).

Alternatively, the second approach employs the Euler-Maruyama (EM) method \cite{Kloeden} and utilizes the weak convergence property \cite{Kloeden} of the EM method. We will adopt the latter approach in our analysis, as it does not require knowledge of advanced stochastic calculus if we accept the weak convergence property. Additionally, this approach can be naturally extended to the analysis on the discrete-time noisy average consensus system.
	Furthermore, the EM method plays a key role in 
	the optimization method to be presented in Section \ref{sec:DU}.
	Our analysis motivates the use EM method for optimizing the covariance.

\subsection{Euler-Maruyama method}

We use the {\em Euler-Maruyama method} 
corresponding to this SDE so as to study
 the stochastic behavior of 
the solution of the SDE (\ref{sde}) defined above.
The EM method is well-known numerical method for 
solving SDEs~\cite{Kloeden}.

Assume that we need  numerical solutions of a SDE
in the time interval $0 \le t \le T$.
We divide this interval into $N$ bins and let
$
t_{k} \equiv k \eta,\ k = 0, 1, \ldots, N
$
where the interval $\eta$ is given by 
$
\eta \equiv {T}/{N}.		
$
Let us define a discretized sample 
$\bm x^{(k)}$ be 
$
\bm x^{(k)} \equiv \bm x(t_k).	
$
It should be noted that, the choice of the width  $\eta$ is crucial 
in order to ensure the stability and the accuracy of the EM method. 
A small width leads to a more accurate solution, but requires more computational time. 
A large width may be computationally efficient but may lead to instability in the solution.

The recursive equation of 
the EM method corresponding to SDE (\ref{sde}) is
given by 
\begin{align}
 \bm x^{(k+1)} = \bm x^{(k)} - \eta \bm L \bm x^{(k)}  + \alpha \bm w^{(k)},\ k = 0,1,2,\ldots, N,
\end{align}
where each element of $\bm w^{(k)} \equiv (w_1^{(k)}, w_2^{(k)},\ldots, w_n^{(k)})^T$ follows
$
	w_i^{(k)} \sim {\cal N}(0, \eta).
$
In the following discussion, we will use the equivalent expression \cite{Kloeden}:
\begin{align} \label{EMscheme}
	\bm x^{(k+1)} = \bm x^{(k)} - \eta \bm L \bm x^{(k)} + \alpha \sqrt{\eta} \bm z^{(k)},\ k = 0,1,2,\ldots, N,
\end{align}
where $\bm z^{(k)}$ is a random vector following the multivariate Gaussian 
distribution ${\cal N}(\bm 0, \bm I)$.
The initial vector $\bm x^{(0)}$ is set to be $\bm c$.
This recursive equation will be referred 
to as the {\em Euler-Maruyama recursive equation}.

\begin{figure}[htbp]
  \centering
  \includegraphics[width=0.92\hsize]{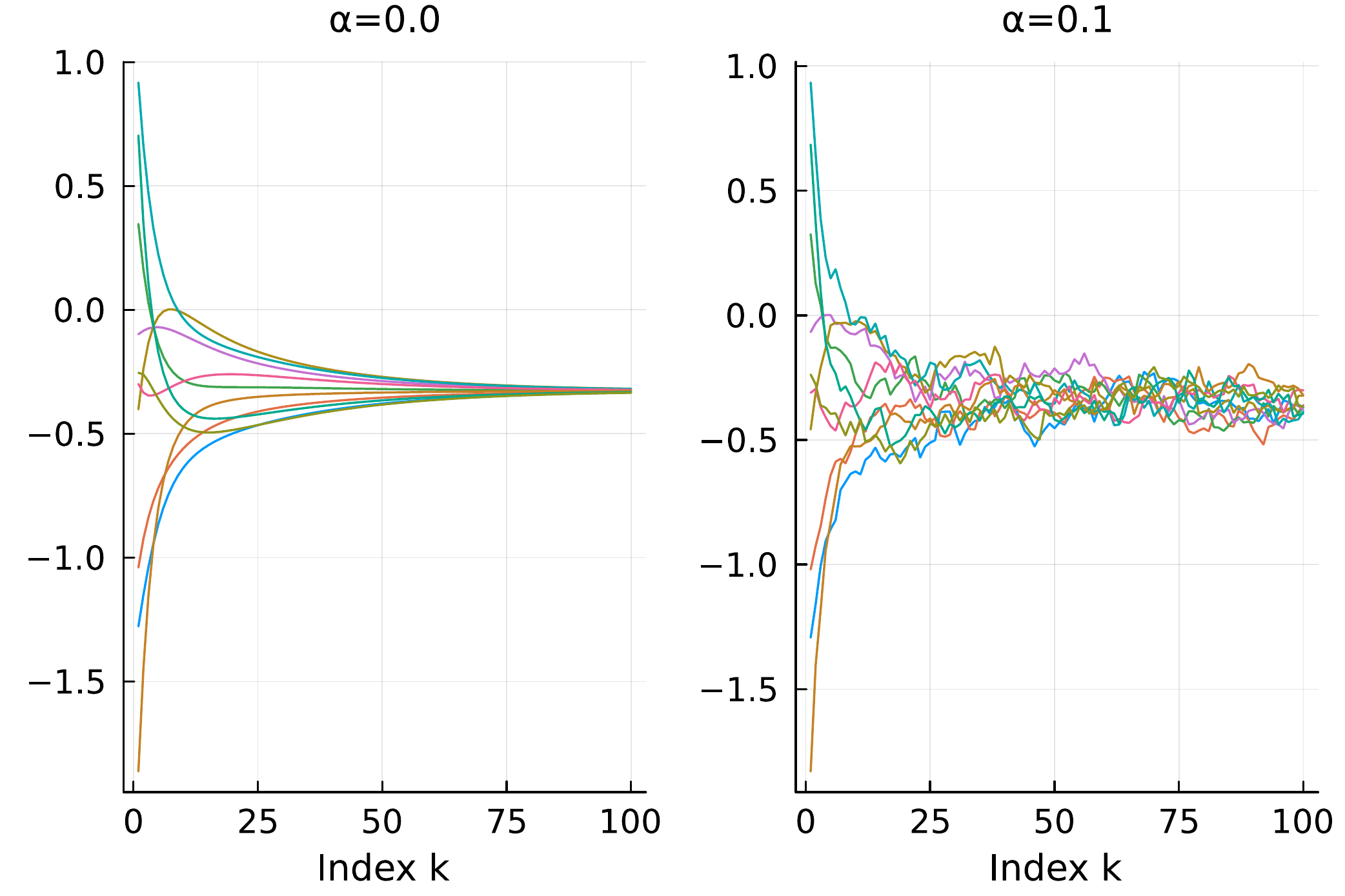}
  \caption{Trajectories of $\bm x(t_k)=(x_1(t_k),\ldots, x_n(t_k))$ estimated by using the EM method. Cycle graph with 10 nodes were used. The range $[0,10.0]$ are discretized with $N=100$ points.
  The consensus average value is $\gamma = -0.3267$. Left panel: noiseless case $(\alpha=0.0)$, Right panel: noisy case $(\alpha=0.1)$.}
  \label{trajectory}
\end{figure}

Figure \ref{trajectory} presents 
a solution evaluated with the EM 
method. The cycle graph with 10 nodes with 
the degree sequence $\bm d  = (2,2,\ldots, 2)$ is assumed.
The initial value is randomly initialized 
as $\bm x(0) \sim {\cal N}(0, \bm I)$.
We can confirm that the state values are certainly converging to the average value $\gamma$
in the case of noiseless case (left). On the other hand, the state vector fluctuates 
around the average in the noisy case (right).

\section{Analysis for Noisy average consensus}
\label{sec:analysis_noisy}

\subsection{Recursive equation for residual error}
In the following, we will analyze the stochastic behavior of the residual error.
This will be the basis for the MSE formula to be presented. 

Recall that the initial state vector 
is $\bm c = (c_1, c_2,\ldots, c_n)^{T}$ and
that the average of the initial values 
is denoted by $\gamma$.
Since the set of eigenvectors $\{\bm \xi_1,\ldots, \bm \xi_n\}$ of $\bm L$ 
is an orthonormal base, we can expand the initial state vector $\bm c$ as 
\begin{align}\label{eigen_extension}
	\bm c = \zeta_1 \bm \xi_1 + \zeta_2 \bm \xi_2 + \cdots + \zeta_n \bm \xi_n,
\end{align}
where the coefficient is obtained by { $\zeta_i = \bm c^{T} \bm \xi_i (i \in [n])$}. 
Note that $\zeta_1 \bm \xi_1 = \gamma \bm 1$ holds.

At the initial index $k=0$,
the Euler-Maruyama recursive equation 
becomes 
\begin{align}
	\bm x^{(1)} = \bm x^{(0)} - \eta \bm L \bm x^{(0)} + \alpha \sqrt{\eta} \bm z^{(0)}.
\end{align}
Substituting (\ref{eigen_extension}) into 
the above equation, we have
\begin{align} \nonumber
	\bm x^{(1)} 
	&=\bm x^{(0)} - \eta \bm L (\zeta_1 \bm \xi_1 + \zeta_2 \bm \xi_2 + \cdots + \zeta_n \bm \xi_n) + \alpha \sqrt{\eta} \bm z^{(0)} \\ \nonumber
	&=\bm x^{(0)} - \eta \zeta_1 \bm L \bm \xi_1 - \eta L (\bm x^{(0)} - \zeta_1\bm \xi_1 ) + \alpha \sqrt{\eta} \bm z^{(0)} \\ 
	&={ \bm x^{(0)}  - \eta \bm L (\bm x^{(0)} - \gamma \bm 1 ) + \alpha \sqrt{\eta} \bm z^{(0)}},
\end{align}
where the equations $L \bm \xi_1 = \bm 0$ and $\zeta_1 \bm \xi_1 = \gamma \bm 1$
are used in the last equality.
Subtracting $\gamma  \bm 1$ from the both sides, we get
\begin{align} \label{initial_cond}
	\bm x^{(1)} - \gamma \bm 1 = (\bm I - \eta \bm L) (\bm x^{(0)} - \gamma \bm 1 )
	+ \alpha \sqrt{\eta} \bm z^{(0)}.
\end{align}
For the index $k \ge 1$, the Euler-Maruyama recursive 
equation can be written as 
\begin{align}
	\bm x^{(k+1)} = (\bm I - \eta \bm L )\bm x^{(k)} + \alpha \sqrt{\eta} \bm z^{(k)}.
\end{align}
{ Subtracting $\gamma  \bm 1$ from the both sides, we have
\begin{align}
	\bm x^{(k+1)} - \gamma \bm 1= (\bm I - \eta \bm L )\bm x^{(k)} - \gamma \bm 1 + \alpha \sqrt{\eta} \bm z^{(k)}.
\end{align}
By using the relation
$
(\bm I - \eta \bm L) \gamma \bm 1= \gamma \bm 1,
$
we can rewrite the above equation as 
\begin{align} \nonumber
	\bm x^{(k+1)} -  \gamma \bm 1 &= (\bm I - \eta \bm L )\bm x^{(k)} - (\bm I - \eta \bm L)  \gamma \bm 1 + \alpha \sqrt{\eta} \bm z^{(k)} \\ \label{k_recursion}
	&= (\bm I - \eta \bm L )(\bm x^{(k)} -  \gamma \bm 1) + \alpha \sqrt{\eta} \bm z^{(k)}.
\end{align}}
It can be confirmed the above recursion (\ref{k_recursion}) is consistent 
with the initial equation (\ref{initial_cond}).
We here summarize the above argument as the following lemma.
\begin{lemma}
	{Let $\bm e^{(k)} \equiv \bm x^{(k)} - \gamma \bm 1$ be 
	the {\em residual error} at index $k$}. 
	The evolution of the residual error of 
	the EM method is described by
\begin{align} \label{error_recursive}
\bm e^{(k+1)}= 
	 (\bm I - \eta \bm L )\bm e^{(k)} + \alpha \sqrt{\eta} \bm z^{(k)}
\end{align}
for $k \ge 0$.
\end{lemma}

The residual error $\bm e^{(k)}$ denotes the error between the average vector $\gamma \bm 1$ and the state vector $\bm x^{(k)}$ at time index $k$. By analyzing the statistical behavior of $\bm e^{(k)}$, we can gain insight into the stochastic properties of the dynamics of the noisy consensus system.

\subsection{Asymptotic mean of residual error}

Let a vector $\bm x \sim {\cal N}(\bm \mu, \bm \Sigma)$.
Recall that the vector obtained by a linear map $\bm y = \bm A \bm x$
also follows the Gaussian distribution, i.e.,  
\begin{align}
\bm y \sim {\cal N}(\bm A\bm \mu, \bm A \bm \Sigma \bm A^{T}),
\end{align}
where $\bm A \in \mathbb{R}^{n \times n}$.
If two Gaussian vectors $\bm a \sim {\cal N}(\bm \mu_a, \bm \Sigma_a)$
and $\bm b \sim {\cal N}(\bm \mu_b, \bm \Sigma_b)$ are independent, 
the sum $\bm z = \bm a + \bm b$ becomes also Gaussian, i.e, 
\begin{align}
\bm z \sim {\cal N}(\bm \mu_a + \bm \mu_b, \bm \Sigma_a + \bm \Sigma_b).
\end{align}

In the recursive equation (\ref{error_recursive}),  
it is evident that $\bm e^{(1)}$ follows a multivariate Gaussian distribution because
\begin{align}
\bm e^{(1)} = 
(\bm I - \eta \bm L )(\bm c - \gamma \bm 1) + \alpha \sqrt{\eta} \bm z^{(0)}
\end{align}
is the sum of a constant vector and a Gaussian random 
vector.
From the above properties of Gaussian random vectors, 
the residual error vector $\bm e^{(k)}$
follows 
the multivariate Gaussian distribution ${\cal N}(\bm \mu^{(k)}, \bm \Sigma^{(k)})$ where the mean vector $\bm \mu^{(k)}$ and
the covariance matrix $\bm \Sigma^{(k)}$ are recursively determined by
\begin{align}
	\bm \mu^{(k+1)} &= (\bm I - \eta  \bm L) \bm \mu^{(k)}, \\ \label{cov_update}
	\bm \Sigma^{(k+1)} &= (\bm I - \eta  \bm L) \bm \Sigma^{(k)}(\bm I - \eta  \bm L)^{T} 
	+ \alpha^2 \eta \bm I
\end{align}
for $k \ge 0$ where 
the initial values are formally given by 
\begin{align}
\bm \mu^{(0)} &= \bm c - \gamma \bm 1, \\
\bm \Sigma^{(0)} &= \bm O.
\end{align}
Solving the recursive equation, we can get 
the asymptotic mean formula as follows.
\begin{lemma}
\label{mean_lemma}
Suppose that $T > 0$ is given. The asymptotic mean 
at $N \rightarrow \infty$ is given by 
\begin{align} \label{asymptotic_mean}
	\lim_{N \rightarrow \infty}\bm \mu^{(N)} 
	= \exp(-\bm L T) (\bm c - \gamma  \bm 1).
\end{align}
\end{lemma}
(Proof) 
The mean recursion is given as
$
	\bm \mu^{(k)} = (\bm I - \eta  \bm L)^k (\bm c - \gamma \bm 1)
$
for $k \ge 1$. Recall that 
the eigenvalue decomposition of $\bm L$
is given by $\bm L = \bm U \mbox{diag}(\lambda_1,\ldots, \lambda_n) \bm U^{T}$. From
\begin{align}
	\bm I - \eta \bm L = \bm U(\bm I - \eta \mbox{diag}(\lambda_1,\ldots, \lambda_n))\bm U^{T},
\end{align}
we have
\begin{align}
	(\bm I - \eta \bm L)^k = \bm U
	\mbox{diag}((1-\eta \lambda_1)^k,\ldots, (1-\eta\lambda_n)^k)
	\bm U^{T}.
\end{align}
This implies, from the definition of exponential function, 
\begin{align}
	\lim_{N \rightarrow \infty} \left(\bm I - \frac{T}{N}\bm L\right)^N = \exp(-\bm L T), 
\end{align}
where $\eta = T/N$. \hfill \fbox{}

It is easy to confirm that the claim of this lemma is consistent with the continuous solution 
of noiseless case (\ref{cont_sol}).
Namely, at the limit of $\alpha \rightarrow 0$, the state evolution of the noisy system 
converges to that of the noiseless system.

\subsection{Asymptotic covariance of residual error}

We here discuss the asymptotic behavior of
the covariance matrix $\bm \Sigma^{(N)}$ at the 
limit of $N \rightarrow \infty$.

\begin{lemma}
\label{cov_lemma}
Suppose that $T > 0$ is given. The asymptotic covariance matrix
at $N \rightarrow \infty$ is given by 
	\begin{align}\label{asymptotic_cov}
		\lim_{N \rightarrow \infty} \bm \Sigma^{(N)}
		= \bm U \mbox{diag}\left(\alpha^2 T, \theta_2,\theta_3, \ldots, \theta_n \right) \bm U^{T},
	\end{align}
where $\theta_i$ is defined by
\begin{align}
\theta_i \equiv \frac{\alpha^2}{2 \lambda_i}\left(1 - e^{-2 \lambda_i T} \right).
\end{align}
	
\end{lemma}
(Proof) Recall that 
\begin{align}
	\bm I - \eta \bm L = \bm U \mbox{diag}(1,1 - \eta \lambda_2\ldots, 1- \eta \lambda_n) \bm U^{T}.	
\end{align}
Let $\bm \Sigma^{(k)} = \bm U \mbox{diag}(s_1^{(k)},\ldots, s_n^{(k)}) \bm U^{T}$.
A spectral representation of the covariance evolution (\ref{cov_update}) is thus given by
\begin{align} \nonumber
	 &\mbox{diag}(s_1^{(k+1)},\ldots, s_n^{(k+1)})  \\
	 &= \mbox{diag}(s_1^{(k)}, s_2^{(k)}(1 - \eta \lambda_2)^2\ldots, s_n^{(k)}(1- \eta \lambda_n)^2)
	 + \alpha^2 \eta \bm I,
\end{align}
where $s_i^{(0)} = 0$.
The first component follows a recursion
$
	s_{1}^{(k+1)} = s_1^{(k)} + \alpha^2\eta 	
$
and thus we have
$
	s_{1}^{(N)} = \alpha^2\eta N = \alpha^2 T.
$
Another component follows
\begin{align}\label{recursive}
	s_i^{(k+1)} = s_i^{(k)}(1 - \eta \lambda_i)^2 + \alpha^2 \eta	.	
\end{align}
Let us consider the characteristic equation of  (\ref{recursive}) which is given by
\begin{align}
	s = s(1 - \eta \lambda_i)^2 + \alpha^2 \eta.
\end{align}
The solution of the equation is given by 
\begin{align}
	s = \frac{\alpha^2 \eta}{1 - (1 - \eta \lambda_i)^2}.	
\end{align}
The above recursive equation (\ref{recursive}) thus can be transformed as
\begin{align}
 s_i^{(k+1)} - s  = (s_i^{(k)} -s )(1 - \eta \lambda_i)^2.	
\end{align}
From the above equation, $s_i^{(N)}$ can be solved as 
\begin{align}
	s_i^{(N)}   = s + (s_i^{(0)} -s )(1 - \eta \lambda_i)^{2N}.	
\end{align}
Taking the limit $N \rightarrow \infty$, we have
\begin{align}
	\lim_{N \rightarrow \infty} s_i^{(N)}   = \frac{\alpha^2}{2 \lambda_i} 
	 \left(1 - e^{-2 \lambda_i T} \right).
\end{align}
We thus have the claim of this lemma.
\hfill \fbox{}


\subsection{Weak convergence of Euler-Maruyama method}

As previously noted,
the asymptotic mean (\ref{asymptotic_mean}) is consistent 
with the continuous solution.
The weak convergence property 
 of the EM method \cite{Kloeden} 
 allows us to obtain the moments of the error at time $t$.

We will briefly explain the weak convergence property.
Suppose a SDE with the form:
\begin{align}
	d \bm x(t) = \phi (\bm x(t))dt + \psi (\bm x(t)) d\bm b(t).
\end{align}
If $\phi$ and $\psi$ are bounded and Lipschitz continuous, then the finite order moment 
estimated by the EM method converges to the exact moment of the solution $\bm x(t)$ 
at the limit $N \rightarrow \infty$ \cite{Kloeden}. This property is called the weak convergence property.
In our case, the SDE (\ref{sde}) has bounded and Lipschitz continuous coefficient functions,
i.e, $\phi(\bm x) =  -\bm L\bm x$ and $\psi(\bm x) = \alpha$. Hence, we can employ the weak convergence property in our analysis.

Suppose $\bm x(t)$ is a solution of SDE (\ref{sde})
with the initial condition $\bm x(0) = \bm c$.
Let $\bm \mu(t)$ be the mean vector of the residual error $\bm e(t) = \bm x(t) - \gamma \bm 1$
and $\bm \Sigma(t)$ is the covariance matrix
of the residual error $\bm e(t)$.

\begin{theorem}
\label{cov_theorem}
For a positive real number $t > 0$, the mean and 
the covariance matrix of { the residual error $\bm e(t)$}
are given by
{
\begin{align}
	\bm \mu(t) &= \exp(-\bm L t) (\bm c - \gamma \bm 1) \\
	\bm \Sigma(t) &= \bm U \mbox{diag}\left(\alpha^2 t, \theta_2, \theta_3, \ldots, \theta_n \right) \bm U^{T}.
\end{align}	}
\end{theorem}
(Proof) Due to the weak convergence property of the
EM method, the first and second moments
of the error are converged to the asymptotic mean 
and covariance of the EM method \cite{Kloeden}, i.e.,
\begin{align}
\bm \mu(T) &= \lim_{N \rightarrow \infty}\bm \mu^{(N)} \\
\bm \Sigma(T) &= \lim_{N \rightarrow \infty}\bm \Sigma^{(N)},
\end{align}
where $N$ and $T$ are related by $T = \eta N$.
Applying Lemmas \ref{mean_lemma} and \ref{cov_lemma} and
replacing the variable $T$ by $t$ 
provide the claim of the theorem. \hfill \fbox{}

\subsection{Mean squared error}
In the following, we assume that 
the initial state vector $\bm c$ follows 
Gaussian distribution ${\cal N}(\bm 0, \bm I)$.

In this setting, $\bm \mu(t)$ also follows 
multivariate Gaussian distribution
with the mean vector $\bm 0$ and 
the covariance matrix $\bm Q(t) \bm Q(t)^{T}$
where 
\begin{align}
	\bm Q(t) \equiv \exp(-\bm L t) \left(\bm I  - \frac 1 n \bm 1 (\bm 1^{T})  \right)
\end{align}
because $\bm \mu(t)$ can be rewritten as
\begin{align}
\bm \mu(t) &= \exp(-\bm L t) (\bm c - \gamma \bm 1) 
= \exp(-\bm L t) \left(\bm I  - \frac 1 n \bm 1 (\bm 1^{T})  \right)\bm c.
\end{align}

\label{MSE_formula}
By using the result of Theorem \ref{cov_theorem},
we immediately have the following corollary indicating the {\em MSE formula}.
\begin{corollary}
The mean squared error (MSE) 
\begin{align} 
{\sf MSE}(t) \equiv {\sf E}[\|\bm x(t) - \gamma \bm 1\|_2^2]	
\end{align}
is given by
\begin{align} \nonumber
{\sf MSE}(t) &=\alpha^2 t + \frac{\alpha^2}{2}\sum_{i=2}^n \frac{1 - e^{-2\lambda_i t}}{\lambda_i} 
+ \mbox{tr}(\bm Q(t) \bm Q(t)^{T}).
\end{align}
\end{corollary}
(Proof)
We can rewrite $\bm x(t)$ as:
\begin{align} \label{channel_model}
\bm x(t) = \gamma \bm 1 + \bm Q(t) \bm c + \bm w,
\end{align}
where $\bm w \sim {\cal N}(\bm 0,  \bm \Sigma(t))$,
and $\bm w$ and $\bm c$ are independent.
We thus have 
\begin{align} \nonumber
{\sf MSE}(t) &= \mbox{tr}(\bm \Sigma(t) ) + \mbox{tr}(\bm Q(t) \bm Q(t)^{T}) \\ \label{VT_formula}
	&=\alpha^2 t + \frac{\alpha^2}{2}\sum_{i=2}^n \frac{1 - e^{-2\lambda_i t}}{\lambda_i} 
+ \mbox{tr}(\bm Q(t) \bm Q(t)^{T})
\end{align}
due to Theorem \ref{cov_theorem}. \hfill \fbox{}

Since the value of the term $\mbox{tr}(\bm Q(t) \bm Q(t)^{T})$ is exponentially decreasing with $t$,
$\mbox{tr}(\bm \Sigma(t) )$ is dominant in ${\sf MSE}(t)$ for sufficiently large $t$. 
For sufficiently large $t$, the MSE is well approximated by the asymptotic MSE (AMSE) as 
\begin{align} \label{AMSE}
{\sf MSE}(t) \simeq {\sf AMSE}(t) \equiv \alpha^2 t + \frac{\alpha^2}{2}\sum_{i=2}^n \frac{1}{\lambda_i}
\end{align}
because 
$\mbox{tr}(\bm Q(t) \bm Q(t)^{T})$ is negligible,
and $1 - e^{-2\lambda_i t}$ can be well approximated to $1$.
We can observe that the sum of inverse eigenvalue $\sum_{i=2}^n ({1}/{\lambda_i})$
of the Laplacian matrix determines the intercept of the ${\sf AMSE}(t)$. In other words, the graph topology 
influences the stochastic error behavior through the sum of inverse eigenvalues of the
Laplacian matrix.

\begin{figure}[htbp]
  \centering
  \includegraphics[width=0.92\hsize]{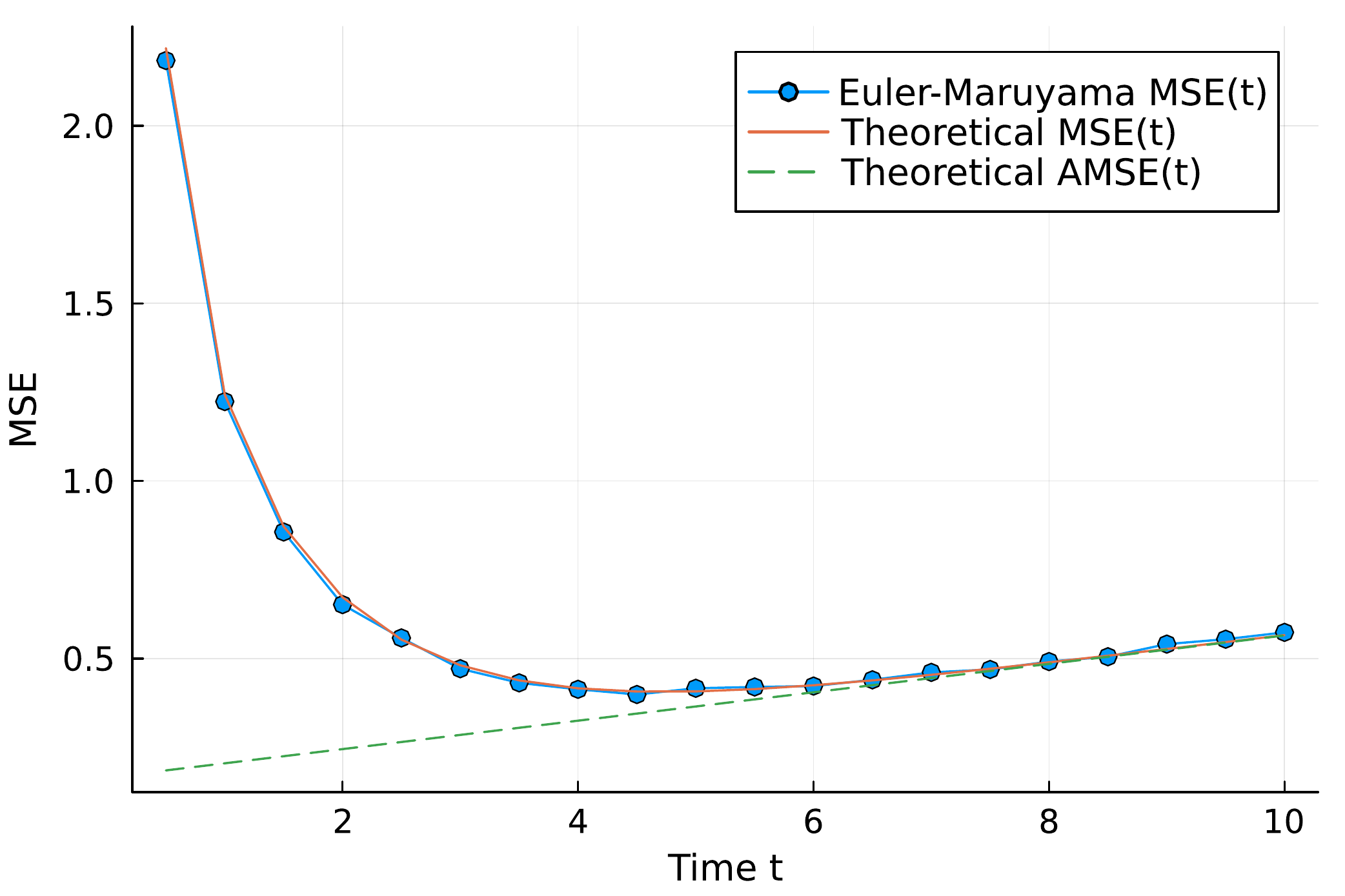}
  \caption{ Comparison of MSE:
  The label Euler-Maruyama represents ${\sf MSE}(t)$ estimated by using samples generated 
  by the EM method.
Theoretical ${\sf MSE}(t)$ represents the values evaluated by (\ref{VT_formula}).Theoretical ${\sf AMSE}(t)$ represents
  the values of ${\sf AMSE}(t)$. 
  Cycle graph with 10 nodes with $\bm d = (2,2,\ldots, 2)$ are used. The parameter setting is as follows: $N = 250$, 
  $T=10$, $\alpha = 0.2$. $5000$ samples are generated by the EM method for estimating ${\sf MSE}(t)$.}
  \label{VandE}
\end{figure}

Figure \ref{VandE} presents a comparison of ${\sf MSE}(t)$ evaluated by
the EM method (\ref{EMscheme}) and the formula in (\ref{VT_formula}).
In this experiment, the cycle graph with 10 nodes is used.
The values of ${\sf AMSE}(t)$ are also included in Fig. \ref{VandE}.
We can see that the theoretical values of ${\sf MSE}(t)$ 
and estimated values by the EM method 
are quite close.

\section{Minimization of mean squared error}
\label{sec:DU}

\subsection{Optimization Problems A and B}


In the previous section, we demonstrated that the MSE 
can be expressed in closed-form. 
It is natural to optimize the edge weights $\{\mu_{ij}\}$ in order to decrease 
the value of the MSE. The optimization of the edge weights is equivalent to 
the optimization of the Laplacian matrix $\bm L$.
There exist several related works that aim to achieve a similar goal for noise-free systems. 
For example, Xiao and Boyd~\cite{Xiao} proposed a method to minimize 
the second eigenvalue to achieve the fastest convergence to the average. 
They formulated the optimization problem 
as a convex optimization problem, which can be solved efficiently.
Kishida et al. \cite{Kishida} presented a deep unfolding-based method 
for optimizing time-dependent edge weights, 
yet these methods are not applicable to systems with noise. 
Optimizing the MSE may be a non-trivial task 
as it involves the sum of the inverse eigenvalues of the Laplacian matrix. 

In this subsection, we will present two optimization problems of edge weights.

\subsubsection{Optimization problem A}

Assume that a degree sequence $\bm d \in \mathbb{R}^+$ is given in advance.
The optimization problem A is the minimization problem of ${\sf MSE}(t^*)$ 
under the given degree sequence where $t^*$ is the predetermined target time given in advance.
The precise formulation of the problem is given as follows:
\begin{align} \nonumber
	\mbox{minimize } {\sf MSE}(t^*) \\ \nonumber
	\mbox{subject to: } \\ 
	\bm L &= \{L_{ij}\} \in \mathbb{R}^{n \times n} \\ \label{symmetric}
	\bm L &= \bm L^T \\ \label{rowsum}
	\bm L &\bm 1 = \bm 0 \\ \label{degreecond}
	&\|\mbox{diag}(\bm L) - \bm d\|_2 < \theta \\ \label{edgecond}
	L_{ij} &= 0,\ (i,j) \notin E.
\end{align}
The constraint (\ref{symmetric}) is imposed for the symmetry of the 
edge weight $\mu_{ij} = \mu_{ji}$ for $(i,j) \in E$. The row sum constraint (\ref{rowsum}) 
is needed for satisfying (\ref{mu_sum}). The constraint (\ref{degreecond}) means that 
$\bm L$ should be close enough to the given degree sequence.
The positive constant $\theta$ can be seen as a tolerance parameter.

One way to interpret the optimization problem A is to consider 
the graph $G$ representing the wireless connection between terminals $i \in [n]$. 
The degree sequence $\bm d = (d_1,d_2,\ldots, d_n)$ can be seen as an allocated receive total wireless power, i.e., 
the terminal $i$ can receive the neighbouring signals up to the total power $d_i$.
If an average consensus protocol is used in such a wireless network 
for specific applications, it is desirable to optimize the ${\sf MSE}(t^*)$ 
while satisfying the power constraint.

\subsubsection{Optimization problem B}

Assume that a real constant $D \in \mathbb{R}^+$ is given in advance.
The optimization problem B is the minimization problem of ${\sf MSE}(t^*)$  
under the situation that the diagonal sum of the Laplacian matrix $\bm L$ is equal to $D$.
The formulation is given as follows:
\begin{align} \nonumber
	\mbox{minimize } {\sf MSE}(t^*) \\ \nonumber
	\mbox{subject to: } \\ 
	\bm L &= \{L_{ij}\} \in \mathbb{R}^{n \times n} \\ 
	\bm L &= \bm L^T \\ 
	\bm L &\bm 1 = \bm 0 \\ \label{diagonal_sum}
	&\left|\sum_{i=1}^n L_{ii} -D \right| < \theta \\
	L_{ij} &= 0,\ (i,j) \notin E.
\end{align}
Following the interpretation above, the power allocation is also optimized in 
this problem.

\subsection{Minimization based on deep-unfolded EM method}

Advances in deep neural networks 
have had a strong impact on the design of 
algorithms for communications and signal processing~\cite{Com1,Com2,Com3}.
{\em Deep unfolding} can be seen as a very effective way
to improve the convergence of iterative algorithms.
Gregor and LeCun introduced the Learned ISTA (LISTA) \cite{LISTA}.
Borgerding et al. also proposed variants of AMP and VAMP 
with trainable capability~\cite{LAMP}\cite{Borgerding}.
Trainable ISTA(TISTA)~\cite{Ito} is another trainable sparse signal recovery algorithm with fast convergence.
TISTA requires only a small number of trainable parameters, which provides a fast and stable training process. 
Another advantage of deep unfolding is that it has a relatively high interpretability of learning results.

The concept behind deep unfolding is rather simple. 
We can embed trainable parameters into 
the original iterative algorithm and then 
unfold the signal-flow graph of the original algorithm. 
The standard supervised training techniques 
used in deep learning, 
such as Stochastic Gradient Descent (SGD) and back propagation, can then be applied
to the unfolded signal-flow graph to optimize the trainable parameters.

The combination of deep unfolding and the differential equation solvers~\cite{Rackauckas}
is a current area of active research in scientific machine learning. 
It should be noted, however, that the technique is not limited to applications within scientific machine learning.
In this subsection, we introduce an optimization algorithm 
that is based on the deep-unfolded EM method. 
The central idea is to use a loss function that approximates ${\sf MSE}(t^*)$.
By using a stochastic gradient descent approach with this loss function, 
we can obtain a near-optimal solution for both optimization problems A and B. 
The proposed method can be easily implemented using any modern neural network 
framework that includes an automatic differentiation mechanism. 
The following subsections will provide a more detailed explanation of the proposed method.

\subsubsection{Mini-batch for optimization}
In an optimization process described below, 
a  number of mini-batches are randomly generated. 
A mini-batch consists of 
\begin{align}
	{\cal M} \equiv \{(\bm c_1, \gamma_1), (\bm c_2, \gamma_2), \ldots 
	(\bm c_K, \gamma_K) \}.
\end{align}
The size parameter $K$ is called the mini-batch size.
The initial value vector $\bm c_i$ follows Gaussian distribution,
i.e., $\bm c_i \sim {\cal N}(\bm 0,\bm I) (i \in [n])$. The corresponding average
value are obtained by $\gamma_i \equiv (1/n)\bm c_i^T \bm I$.

\subsubsection{Loss function for Optimization problem A}
The loss function corresponding to a mini-batch ${\cal M}$ is given by 
\begin{align}
	E_{\cal M}(\bm L) \equiv \frac 1 K \sum_{i = 1}^K \| \bm \chi(\bm c_i) - \gamma_i \bm 1\|^2_2 + P_A(\bm L),
\end{align}
where $\bm \chi(\bm c_i) \equiv \bm x^{(N)}$ is the random variable given by
the Euler-Maruyama recursion:
\begin{align} 
	\bm x^{(k+1)} = \bm x^{(k)} - \eta \bm L \bm x^{(k)} + \alpha \sqrt{\eta} \bm z^{(k)},\ k = 0,1,2,\ldots, N,
\end{align}
with $\bm x^{(0)} = \bm c_i$.
The first term of the loss function can be regarded as an 
approximation of ${\sf MSE}(t^*)$:
\begin{align}
\frac 1 K \sum_{i = 1}^K \| \bm \chi(\bm c_i) - \gamma_i \bm 1\|^2_2 \simeq {\sf MSE}(t^*)
\end{align}
for sufficiently large $K$ and $T = t^*$.

The function $P_A(\bm L)$ is a penalty function corresponding to the
constraints (\ref{symmetric})--(\ref{edgecond}) defined by
\begin{align} \nonumber
P_A(\bm L) &\equiv \rho_1 \|\bm L - \bm L^T\|_F^2 + \rho_2\|\bm L \bm 1\|_2^2 
	+ \rho_3 \| \mbox{diag}(\bm L) - \bm d \|_2^2 \\
	+ \rho_4 \|\bm L 	&\odot \bm M\|_F^2,
\end{align}
where $\bm M = \{M_{ij}\}$ is the mask matrix defined by
\begin{align}
	M_{ij} \equiv 
	\left\{
	\begin{array}{cc}
		1, & (i,j) \notin E \\
		0, & \mbox{otherwise}.
	\end{array}
	\right.
\end{align}
The operator $\odot$ represents the Hadamard matrix product.
The positive constants $\rho_i (i \in [4])$ controls relative strength of each penalty term.
The first term of the penalty function corresponds to the 
symmetric constraint (\ref{symmetric}). The term $\|\bm L \bm 1\|_2^2$
is the penalty term for the row sum constraint (\ref{rowsum}).
The third term $\| \mbox{diag}(\bm L) - \bm d \|_2^2$ is included for the degree constraint.
The last term $\|\bm L \odot \bm M\|_F^2$ enforces $L_{ij}$ to be very small
if $(i,j) \notin E$.

Due to these penalty terms in $P_A(\bm L)$,  
the violations on the constraints (\ref{symmetric})--(\ref{edgecond}) are suppressed 
in an optimization process.

\subsubsection{Loss function for Optimization problem B}

For Optimization problem B, we use almost the same same loss function:
\begin{align}
	E_{\cal M}(\bm L) \equiv \frac 1 K \sum_{i = 1}^K \| \bm \chi(\bm c_i) - \gamma_i \bm 1\|^2_2 + P_B(\bm L).
\end{align}
In this case, we use the penalty function matched to the feasible conditions of Optimization problem B:
\begin{align} \nonumber
P_B(\bm L) &\equiv \rho_1 \|\bm L - \bm L^T\|_F^2 + \rho_2 \|\bm L \bm 1\|_2^2 
	+ \rho_3 \left(\sum_{i=1}^n L_{ii} - D \right)^2 \\
	+ \rho_4 \|\bm L 	&\odot \bm M\|_F^2.
\end{align}
The third term of $P_B(\bm L)$ corresponds to the diagonal sum condition (\ref{diagonal_sum}).

\subsubsection{Optimization process}

\begin{figure*}[htbp]
  \centering
  \includegraphics[width=0.8\hsize]{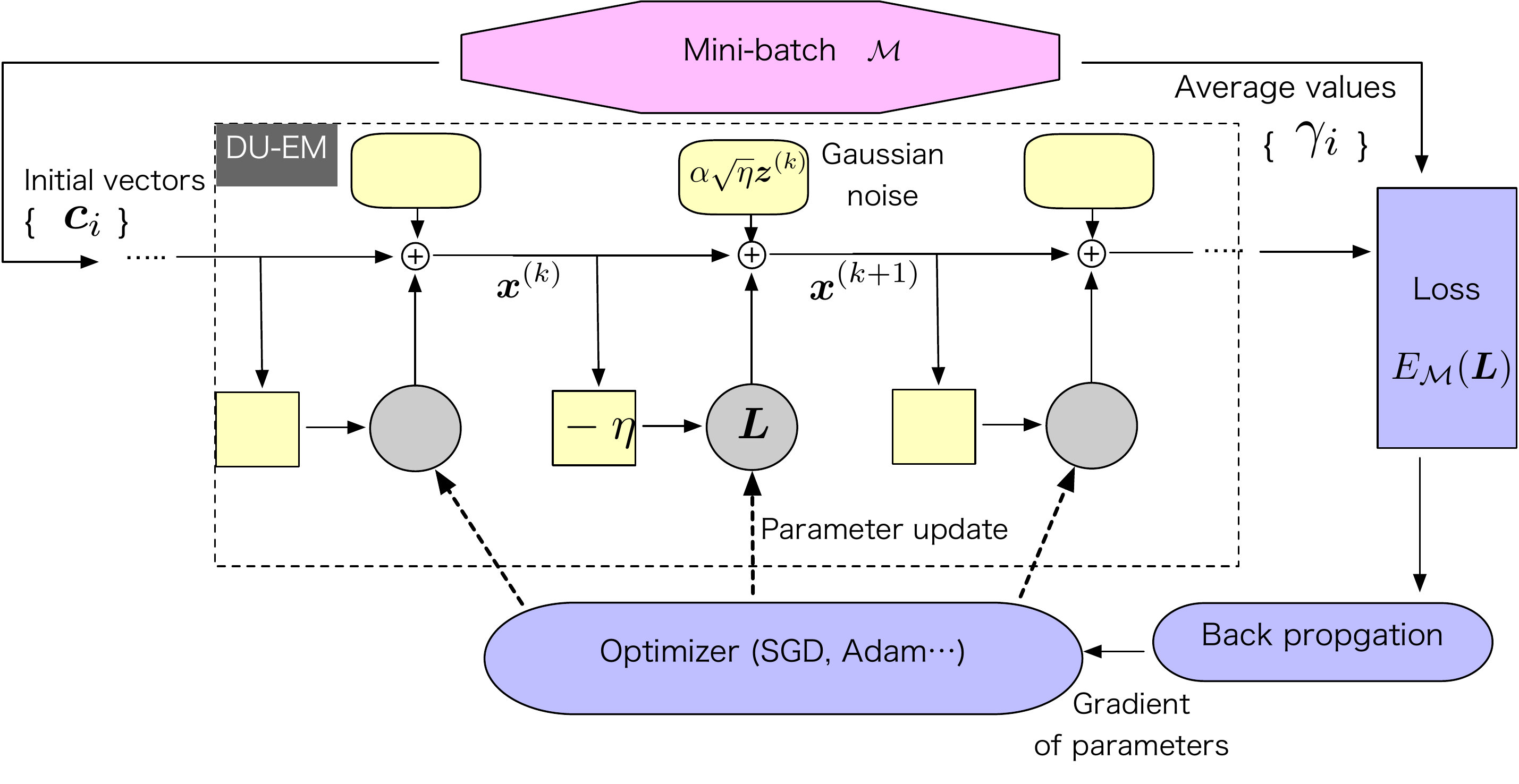}
  \caption{Block diagram of optimization process in Algorithm \ref{opt_alg}. The core of the algorithm is
  the DU-EM method for approximating ${\sf MSE}(t^*)$. Several standard deep learning techniques such as back propagation and stochastic gradient descent can be applied to update the trainable matrix $\bm L$.}
  \label{DU_EM}
\end{figure*}

The optimization process is summarized in Algorithm \ref{opt_alg}.
This optimization algorithm is mainly based on 
the Deep-unfolded Euler-Maruyama (DU-EM) method for approximating ${\sf MSE}(t^*)$.
The initial value of the matrix $\bm L$ is assumed to be
the $n \times n$ zero matrix $\bm O^{n \times n}$.
The main loop can be regarded as a stochastic gradient descent 
method minimizing the loss values. The update of $\bm L$ (line 5) can 
be done by any optimizer such as the Adam optimizer.
The gradient of the loss function (line 4) can be easily evaluated by using 
an automatic differentiation mechanism included in recent neural network frameworks 
such as TensorFlow, PyTorch, Jax, and Flux.jl with Julia.
The block diagram of the Algorithm \ref{opt_alg} is shown in Fig. \ref{DU_EM}.
\begin{algorithm}
 \caption{Optimization process using DU-EM method}
 \label{opt_alg}
 \begin{algorithmic}[1]
 \renewcommand{\algorithmicrequire}{\textbf{Input:}}
 \renewcommand{\algorithmicensure}{\textbf{Output:}}
  \REQUIRE graph $G$, tolerance $\theta$, degree sequence $\bm d$ or degree sum $D$
  \ENSURE Laplacian matrix $\bm L_{out}$
  \STATE Let $\bm L \equiv \bm O^{n \times n}$
   \FOR {$i = 1$ to $I$}
  	\STATE Generate a mini-batch ${\cal M}$ randomly.
	\STATE Compute the gradient of the loss function
	\begin{align}\nonumber
	\bm g \equiv \nabla E_{\cal M}(\bm L)
	\end{align}
	\STATE The matrix $\bm L$ is updated by using $\bm g$.
  \ENDFOR
  \STATE $\bm L_{out} \equiv \mbox{round}_{\theta, *}(\bm L)$
 \end{algorithmic} 
 \end{algorithm}

The stochastic optimization process outlined in Algorithm \ref{opt_alg} 
is unable to guarantee that the obtained solution will be strictly feasible. 
To ensure feasibility, it is necessary to search for a feasible solution 
that is near the result obtained by optimization. 
This is accomplished by using the {\em round function} 
$\mbox{round}_{\theta, *}(\cdot)$ at line 7 of Algorithm \ref{opt_alg}.

The specific details for the round function used for optimization problem A are outlined in Algorithm \ref{round_A}. The first step in the algorithm, $\bm L \equiv (\bm L_{in} + \bm L_{in}^T)/2$, ensures that the resulting matrix is symmetric. 
The nested loop from line 2 to line 7 is used to enforce the degree constraint and the constraint $L_{ij} = 0\  (i,j) \notin E$. The single loop from line 9 to line 11 is implemented to satisfy the constraint $\bm L \bm 1 = \bm 0$. The output of the round function $\mbox{round}_{\theta, \bm d}(\cdot)$ guarantees that the constraints (\ref{symmetric})-(\ref{edgecond}) of optimization problem A are strictly satisfied. A similar round function can be constructed for optimization problem B, which is presented in Algorithm \ref{round_B}.

\begin{algorithm}
 \caption{Round function $\mbox{round}_{\theta, \bm d}(\cdot)$ for Opt. prob. A}
 \label{round_A}
 \begin{algorithmic}[1]
 \renewcommand{\algorithmicrequire}{\textbf{Input:}}
 \renewcommand{\algorithmicensure}{\textbf{Output:}}
  \REQUIRE Matrix $\bm L_{in}$, degree sequence $\bm d$, threshold value $\theta$
  \ENSURE  Laplacian matrix $\bm L_{out}$ satisfying (\ref{symmetric})--(\ref{edgecond})
  \STATE Let $\bm L \equiv (\bm L_{in} + \bm L_{in}^T)/2$
   \FOR {$i = 1$ to $n$}
  	\STATE $L_{ii} \equiv d_i$
  	\FOR{$j = 1$ to $n$}
  		\STATE If $(i,j) \notin E$, then let $L_{ij} \equiv 0$
  	\ENDFOR
  \ENDFOR
  \STATE $\bm \epsilon = (\epsilon_1,\ldots, \epsilon_n)^T \equiv \bm L \bm 1$
  \FOR {$i = 1$ to $n$}
   \STATE Let $L_{ii} \equiv L_{ii} - \epsilon_i$ 
  \ENDFOR
  \IF{$\|\mbox{diag}(\bm L) - \bm d\|_2 \ge \theta$}
  \STATE Quit with declaration ``optimization failed''
  \ENDIF
  \STATE Output $\bm L_{out} \equiv \bm L$
 \end{algorithmic} 
 \end{algorithm}

 \begin{algorithm}
 \caption{Round function $\mbox{round}_{\theta, D}(\cdot)$ for Opt. prob. B}
 \label{round_B}
 \begin{algorithmic}[1]
 \renewcommand{\algorithmicrequire}{\textbf{Input:}}
 \renewcommand{\algorithmicensure}{\textbf{Output:}}
  \REQUIRE Matrix $\bm L_{in}$, degree sum $D$, threshold value $\theta$
  \ENSURE  Laplacian matrix $\bm L_{out}$ 
  \STATE Let $\bm L \equiv (\bm L_{in} + \bm L_{in}^T)/2$
   \FOR {$i = 1$ to $n$}
  	\FOR{$j = 1$ to $n$}
  		\STATE If $(i,j) \notin E$, then let $L_{ij} \equiv 0$
  	\ENDFOR
  \ENDFOR
  \STATE $\bm \epsilon = (\epsilon_1,\ldots, \epsilon_n)^T \equiv \bm L \bm 1$
  \FOR {$i = 1$ to $n$}
   \STATE Let $L_{ii} \equiv L_{ii} - \epsilon_i$ 
  \ENDFOR
  \IF{$\left|\sum_{i=1}^n L_{ii} -D \right| \ge \theta$}
  \STATE Quit with declaration ``optimization failed''
  \ENDIF
  \STATE Output $\bm L_{out} \equiv \bm L$
 \end{algorithmic} 
 \end{algorithm}

 \section{Numerical results} 

\subsection{Choice of Number of bins for EM-method}

In the previous sections, we proposed a DU-based optimization method.
This section presents results of numerical experiments.
For these experiments, we used the automatic differentiation mechanism provided by
{Flux.jl} \cite{Flux} on Julia Language \cite{Julia}.

Before discussing the optimization of MSE, 
we first examine the choice number of bins, $N$.
Small $N$ is beneficial for computational efficiency but
it may lead to inaccurate estimation of MSE.
In this subsection, we will compare the Monte carlo estimates of MSE 
estimated by the EM-method.

The Karate graph is a well-known graph of a small social network. 
It represents the relationships between 34 members of a karate club 
at a university. The graph consists of 34 nodes, 
which represent the members of the club, and 78 edges,
which represent the relationships between the members.

\begin{figure}[htbp]
  \centering
  \includegraphics[width=0.92\hsize]{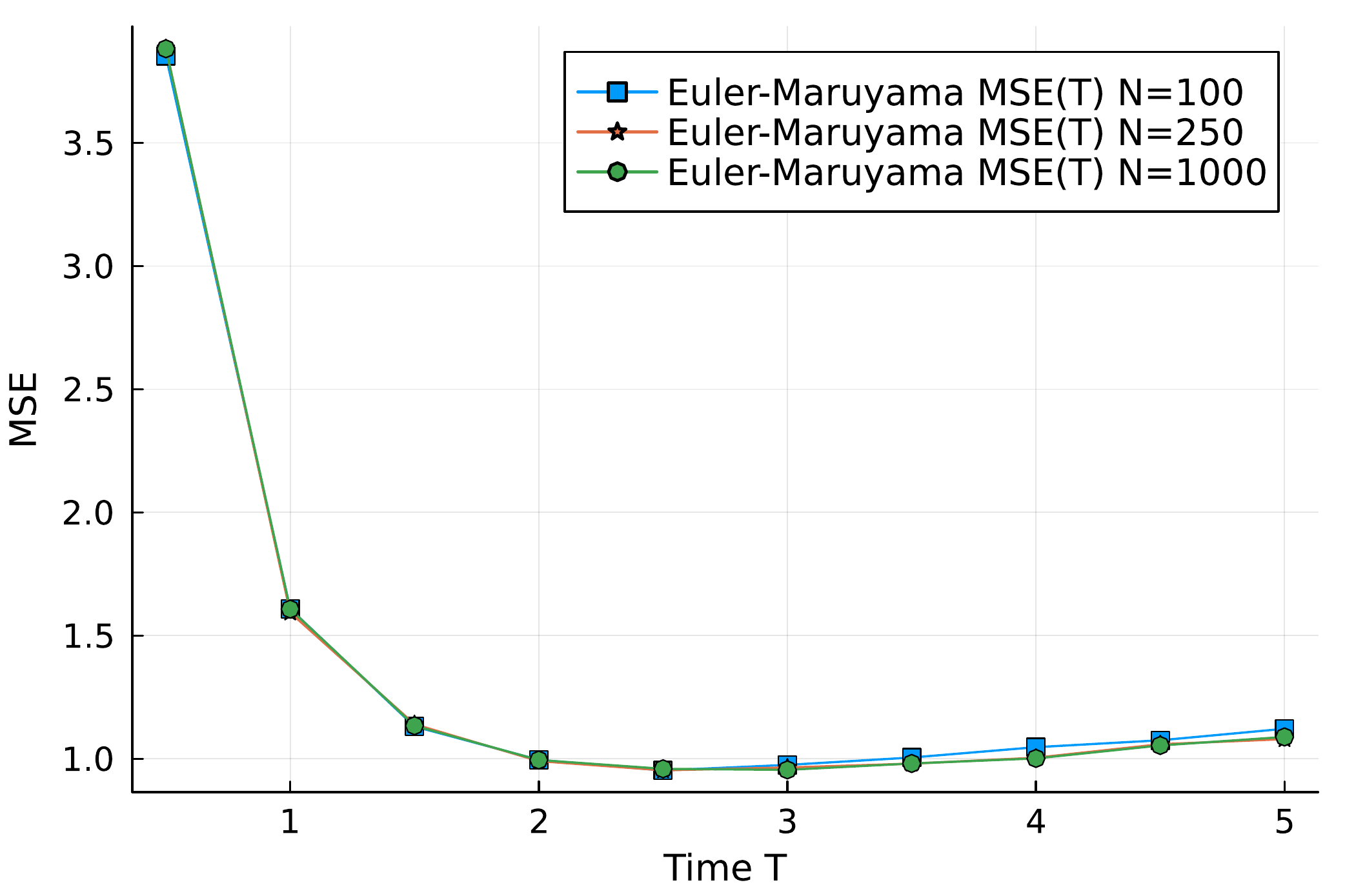}
  \caption{MSE values estimated by Monte Carlo method based on the EM method. Karate graph ($n=34$) and its unweighted Laplacian is assumed.}
  \label{N_dependency}
\end{figure}

Figure \ref{N_dependency} compares three cases, i.e., $N = 100, 250, 1000$.
No visible difference can be observed in the range from $T=0$ to $T=5$.
In the following experiments, we will use $N = 250$ for EM-method.

\subsection{Petersen graph  (Optimization problem A)}

\begin{figure}[htbp]
  \centering
  \includegraphics[width=0.92\hsize]{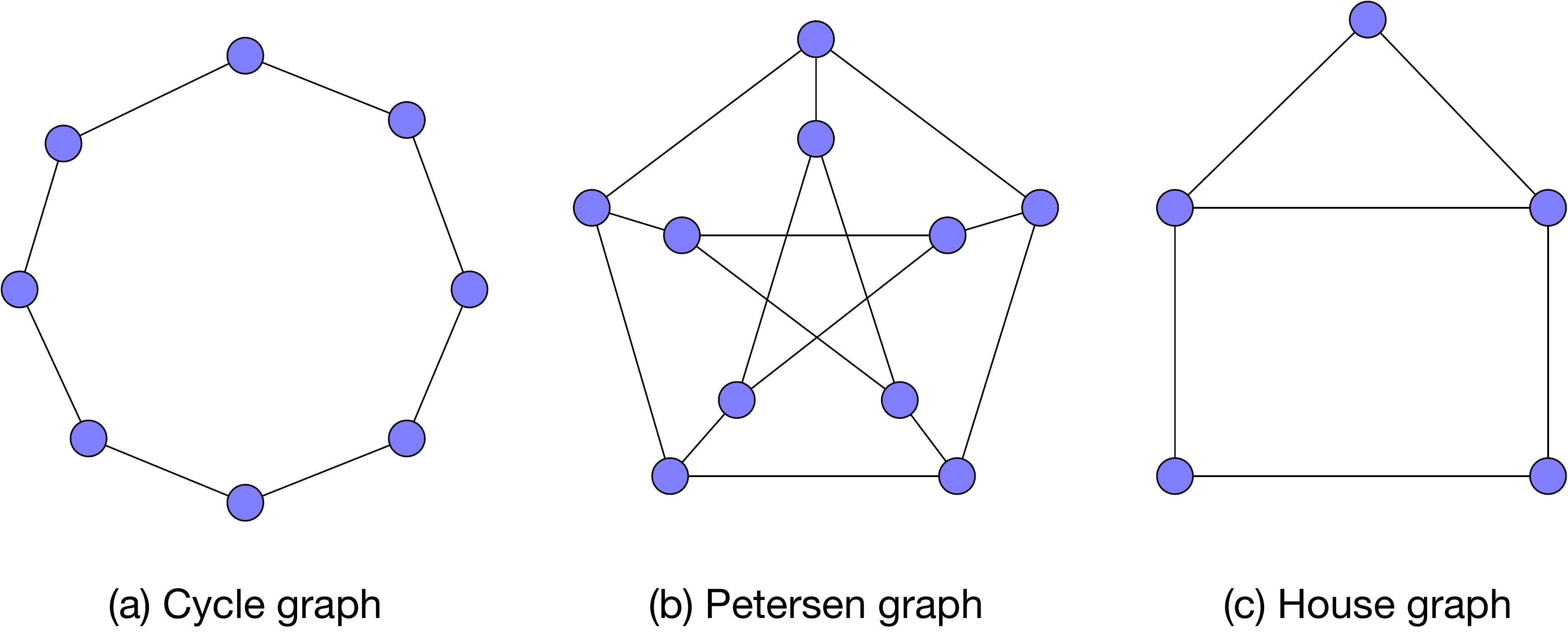}
  \caption{Small graphs: (a) Cycle graph, (b) Petersen graph, (c) House graph.}
  \label{small_graphs}
\end{figure}

Petersen graph is a 3-regular graph with $n=10$ nodes (Fig.\ref{small_graphs}(a)).
In this subsection, we will examine the behavior of our optimization 
algorithm of ${\sf MSE}(t)$ for Petersen graph.

An adjacency matrix $\bm A \equiv \{A_{ij}\} \in \mathbb{R}^{n \times n}$ of a graph $G \equiv (V,E)$ is
defined by
\begin{align}
	A_{ij} \equiv 
	\left\{
	\begin{array}{cc}
		1, & (i,j) \in E \\
		0, & \mbox{otherwise}.
	\end{array}
	\right.
\end{align}
An unweighted Laplacian matrix $\bm L$ is defined by
\begin{align}
\bm L \equiv \bm D - \bm A,
\end{align}
The degree matrix $\bm D = \{D_{ij}\}$ is a diagonal matrix where $D_{ii}$ is the degree of the node $i$.
Namely, an unweighted Laplacian corresponds the case where $\mu_{ij} = \mu_{ji} = 1$ for any $(i,j) \in E$.

In the following discussion, let $\bm L_{P}$ be the unweighted Laplacian matrix of Petersen graph.
We assume Optimization problem A with the degree sequence $\bm  d \equiv \mbox{diag}(\bm L_{P}) = (3,3,\ldots, 3)$.

The parameter setting is as follows.
The mini-batch size is set to $K = 25$. The noise intensity is $\alpha = 0.3$.
The penalty coefficients are $\rho_1 = \rho_2 = \rho_3 = \rho_4 = 10$.
For time discretization, we use $T = 4, N = 250$. The number of iterations for an optimization process 
is set to 3000. The tolerance parameter is set to $\theta = 0.1$.
In the optimization process, we used the Adam optimizer with a learning rate of 0.01.

\begin{figure}[htbp]
  \centering
  \includegraphics[width=0.92\hsize]{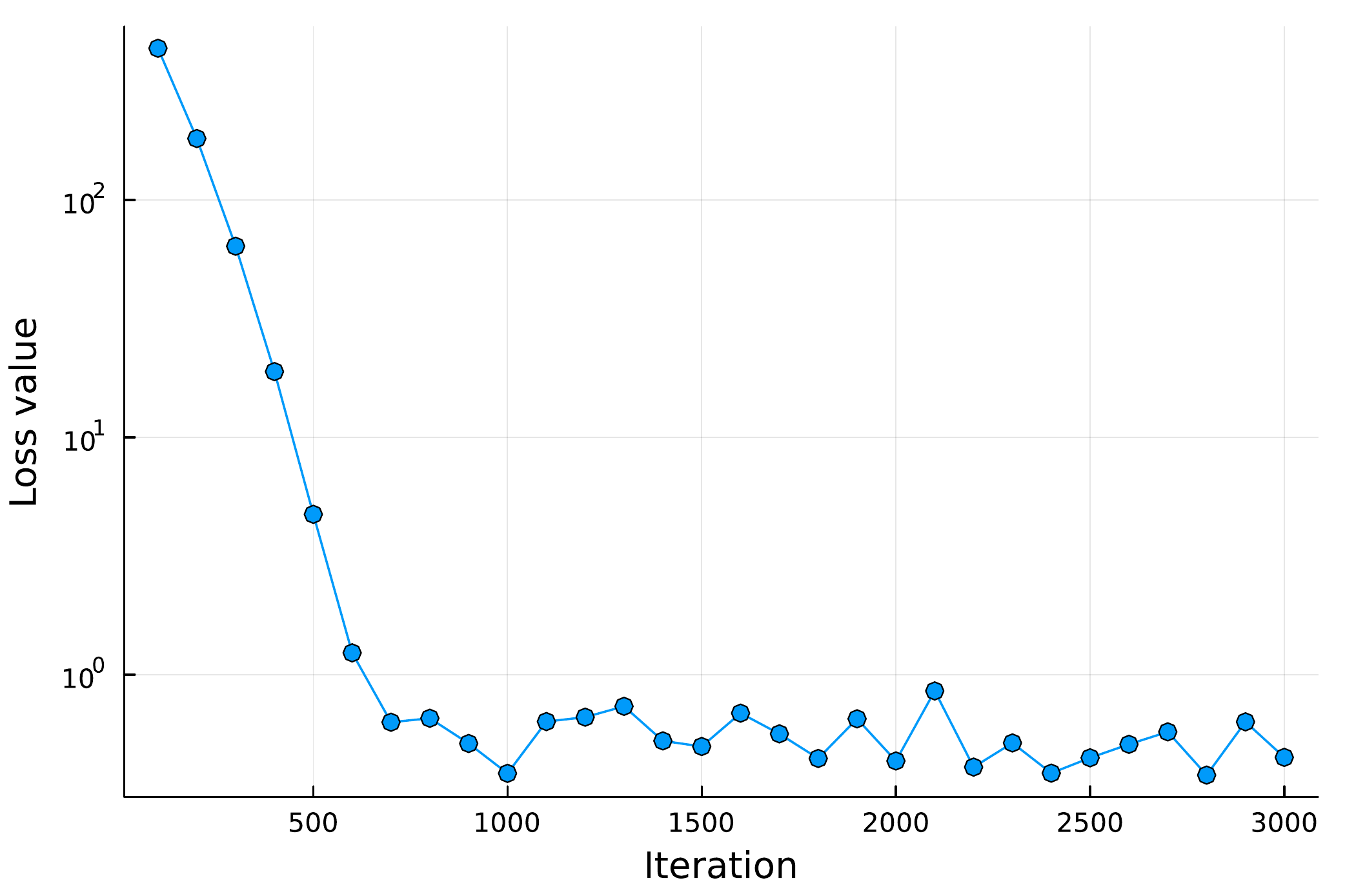}
  \caption{Loss values  in an optimization process: Optimization prob. A for Petersen graph}
  \label{Petersen_loss}
\end{figure}

The loss values of an optimization process of Algorithm \ref{opt_alg} are presented 
in Fig.\ref{Petersen_loss}. In the initial stages of the optimization process, 
the loss value is relatively high since the initial $\bm L$ is set to the zero matrix,
which means that the system cannot achieve average consensus.
The loss value decreases monotonically until around iteration 700, after which it fluctuates within a range of $700 \le k \le 3000$. The graph shows that the matrix $\bm L$ 
in Algorithm \ref{opt_alg} is being updated appropriately and that the loss value, which approximates ${\sf MSE}(t)$, is decreasing.

Let us denote the Laplacian matrix obtained by the optimization process as $\bm L^*$.
Table \ref{quantities_Lstar} summarizes several important quantities regarding $\bm L^*$.
The top 4 rows of Table \ref{quantities_Lstar} indicate that $\bm L^*$ is certainly a feasible 
solution satisfying (\ref{symmetric})--(\ref{edgecond}) because we set $\theta = 0.1$.
This numerical results confirms that the round function $\mbox{round}_{\theta, \bm d}(\cdot)$ works appropriately.
The last row of Table \ref{quantities_Lstar} shows that $\bm L^*$ is very close to the unweighted Laplacian
matrix $\bm L_P$. Since Petersen graph is regular and has high symmetry, it is conjectured that
$\bm L_P$ is the optimal solution for Optimization problem A. Thus, the closeness between $\bm L_P$ and
$\bm L^*$ can be seen as a convincing result.

\begin{table}[htbp]
\centering
\caption{Several quantities on optimization result $\bm L^*$}
\label{quantities_Lstar}
\begin{tabular}{cc}
\hline
\hline
$\|\bm L^* - \bm L^{*T} \|_F$ & 0 \\
$\|\bm L^* \bm 1 \|_2$ & 0 \\
$\|\bm L^* \odot M \|_F$ & 0 \\
$\|\mbox{diag}(\bm L^*) - \bm d \|_2$ & $5.74 \times 10^{-3}$ \\
$\|\bm L_P - \bm L^*\|_F$ & 0.188 \\
\hline 
\end{tabular}
\end{table}

\begin{figure}[htbp]
  \centering
  \includegraphics[width=0.92\hsize]{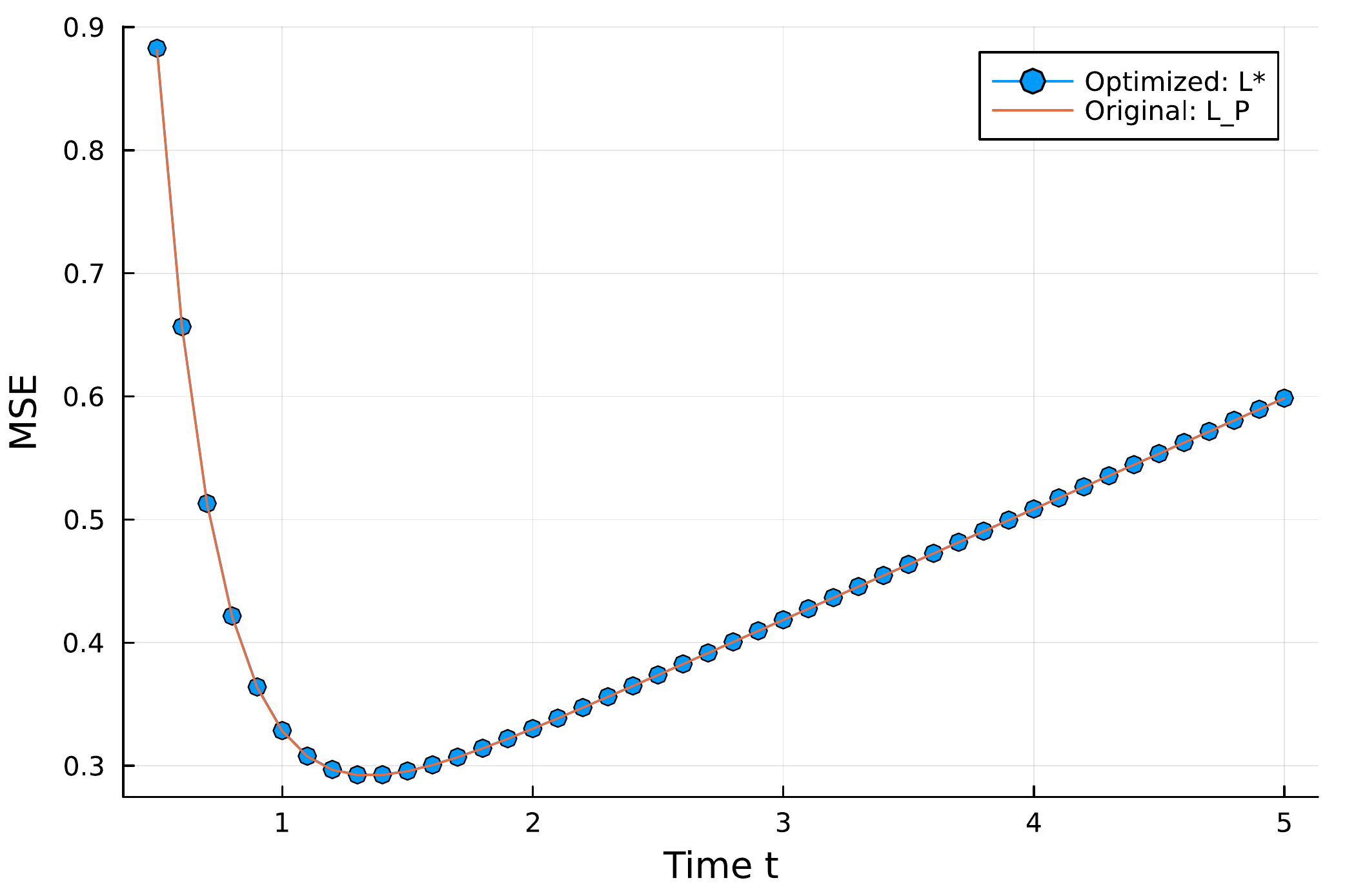}
  \caption{Petersen graph: MSE values of the optimization result $\bm L^*$ and the unweighted Laplacian matrix $\bm L_P$.}
  \label{MSE_Petersen}
\end{figure}

The MSE values of the optimization result $\bm L^*$ and the unweighted Laplacian matrix $\bm L_P$
are presented in Fig.\ref{MSE_Petersen}. These values are evaluated by the MSE formula (\ref{VT_formula}).
No visible difference can be seen between two curves. This means that Algorithm \ref{opt_alg} 
successfully found a good solution for Optimization problem A in this case.

\subsection{Karate graph (Optimization problem A)}

We here consider Optimization problem A on the Karate graph.
Let $\bm L_K$ be the unweighted Laplacian matrix of the Karate graph.
The target degree sequence is set to 
$
\bm d \equiv \mbox{diag}(\bm L) =(16, 9, 10,\ldots,12, 17).
$
The parameter setting for an optimization process is given as follows.
The mini-batch size is set to $K = 50$. 
The noise intensity is set to $\alpha = 0.3$.
The penalty coefficients are $\rho_1 = \rho_2 = \rho_3 = \rho_4 = 10$.
We use $T = 2, N = 250$ for DU-EM method.
The number of iterations for an optimization process 
is set to 5000. The tolerance is set to $\theta = 0.1$.
In the optimization process, 
we used the Adam optimizer with learning rate 0.01.

\begin{figure}[htbp]
  \centering
  \includegraphics[width=0.95\hsize]{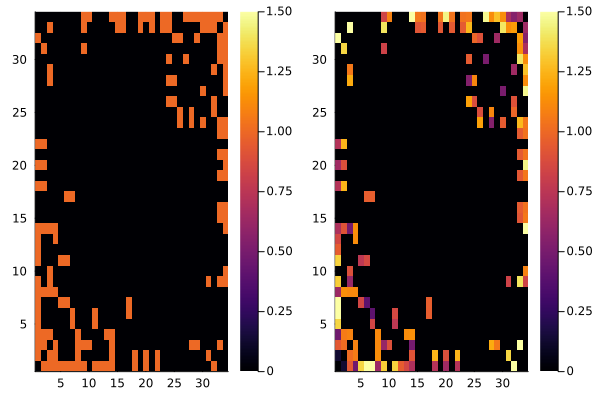}
  \caption{Absolute values of non-diagonal elements in $\bm L_K$ and $\bm L^*$: (left panel) Laplacian matrix ${\bm L}_K$ of the Karate graph, (right panel) The Laplacian matrix $\bm L^*$ obtained by an optimization process.}
  \label{Karate_heatmap}
\end{figure}

Assume that $\bm L^*$ is the Laplacian matrix obtained by
an optimization process. The matrix $\bm L^*$ is a feasible solution satisfying 
all the constraints (\ref{symmetric})--(\ref{edgecond}).
For example, we have $\|\mbox{diag}(\bm L^*) - \bm d\|_2 = 0.0894 < 0.1$.
Figure \ref{Karate_heatmap} presents the absolute values of non-diagonal elements 
in $\bm L_K$ and $\bm L^*$. According to its definition, 
the absolute value of a non-diagonal element of $\bm L_K$ take the value one (left panel). On the other hand,
we can observe that non-diagonal elements of $\bm L^*$ takes the absolute values
in the range $0$ to $1.5$.

\begin{figure}[htbp]
  \centering
  \includegraphics[width=0.95\hsize]{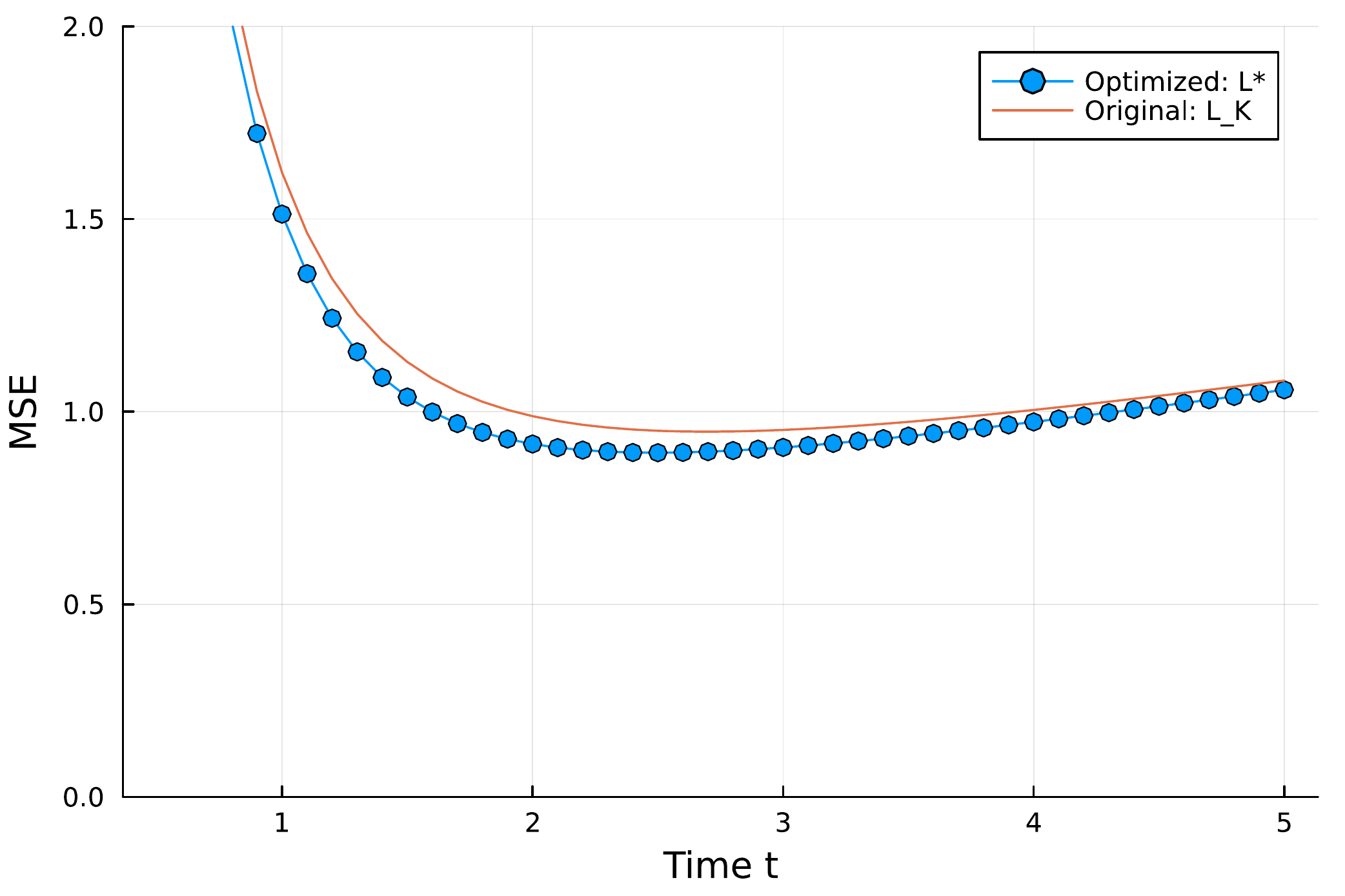}
  \caption{Karate graph: MSE values of $\bm L^*$ ($\bm d = \mbox{diag}(\bm L_K)$) and the unweighted Laplacian matrix $\bm L_K$.}
  \label{MSE_Karate}
\end{figure}

We present the MSE values of the optimization result $\bm L^*$ and the unweighted Laplacian matrix $\bm L_K$ in Fig.\ref{MSE_Karate}. 
These values are evaluated by the MSE formula (\ref{VT_formula}).
It can be seen that the optimized Laplacian $\bm L^*$ provides 
smaller MSE values. In this case, appropriate assignment of weights $\mu_{ij}$
improves the noise immunity of the system.
The inverse eigenvalue sums of the Laplacian matrices 
$\bm L_K$ and $\bm L^*$ are 13.83 and 13.41, respectively.
In this case, the optimization process of Algorithm \ref{opt_alg} 
can successfully provide a feasible Laplacian matrix 
with smaller inverse eigenvalue sum.
As shown in (\ref{AMSE}), the inverse eigenvalue sum determines 
the behavior of ${\sf MSE}(t)$. 


\subsection{House graph (Optimization problem B)}

The house graph (Fig.\ref{small_graphs}(c)) is a small 
irregular graph with 5 nodes defined by
the adjacency matrix:
\begin{align}
\bm A = 
\begin{pmatrix}
0 & 1 & 1 & 0 & 0 \\
1 & 0 & 0 & 1 & 0 \\
1 & 0 & 0 & 1 & 1 \\
0 & 1 & 1 & 0 & 1 \\
0 & 0 & 1 & 1 & 0 \\
\end{pmatrix}.
\end{align}
We thus have the unweighted Laplacian $\bm L_H$ of the house graph as
\begin{align}
\bm L_H =
\begin{pmatrix}
2 & -1 & -1 & 0 & 0 \\
-1 & 2 & 0 & -1 & 0 \\
-1 & 0 & 3 & -1 & -1 \\
0 & -1 & -1 & 3 & -1 \\
0 & 0 & -1 & -1 & 2 \\
\end{pmatrix},
\end{align}
where the diagonal sum of $\bm L_H$ is 12.

We made two optimizations for $D = 12$ and $D = 20$.
The parameter setting is almost the same as the one used in the previous 
subsection. Only the difference is to use $\rho_3 = 0.1$ as the
diagonal sum penalty constant.
As results of the optimization processes, 
we have two Laplacian matrices $\bm L^*_{12}$ $(D=12)$ and
$\bm L^*_{24}$ $(D=24)$ as follows:
\begin{align}
\bm L^*_{12} = 
\begin{pmatrix}
 2.29  &-1.05 &  -1.23 &  0   & 0 \\
 -1.05 & 2.29 & 0      & -1.24& 0 \\
 -1.23 & 0    & 2.70   &-0.44&-1.03 \\
 0     & -1.24&  -0.44&   2.71&   -1.03 \\
0     & 0    &  -1.03 &  -1.03 &    2.06 \\
\end{pmatrix}
\end{align}
\begin{align}
\bm L^*_{24} = 
\begin{pmatrix}
  4.80 &  -2.70 & -2.09 &   0 &    0 \\
 -2.70 &  4.79  & 0        & -2.08&   0 \\
 -2.09 &  0        &  4.81  &  -0.37  & -2.35 \\
  0       &  -2.08 & -0.37&   4.85 &  -2.40 \\
  0       &   0       & -2.35 &  -2.40 &   4.74\\
\end{pmatrix}
\end{align}
The diagonal sums of $\bm L^*_{12}$ and $\bm L^*_{24}$ are 12.04 and 23.99, respectively. 
Compared with $\bm L^*_{12}$ with $\bm L_H$, the diagonal elements of 
$\bm L^*_{12}$ are more flat:
\begin{align}
	\mbox{diag}(\bm L^*_{12}) &= (2.29,2.29,2.70, 2.71,2.06)^T, \\
	\mbox{diag}(\bm L_{H}) &= (2,2,3, 3,2)^T. 
\end{align}

\begin{figure}[htbp]
  \centering
  \includegraphics[width=0.95\hsize]{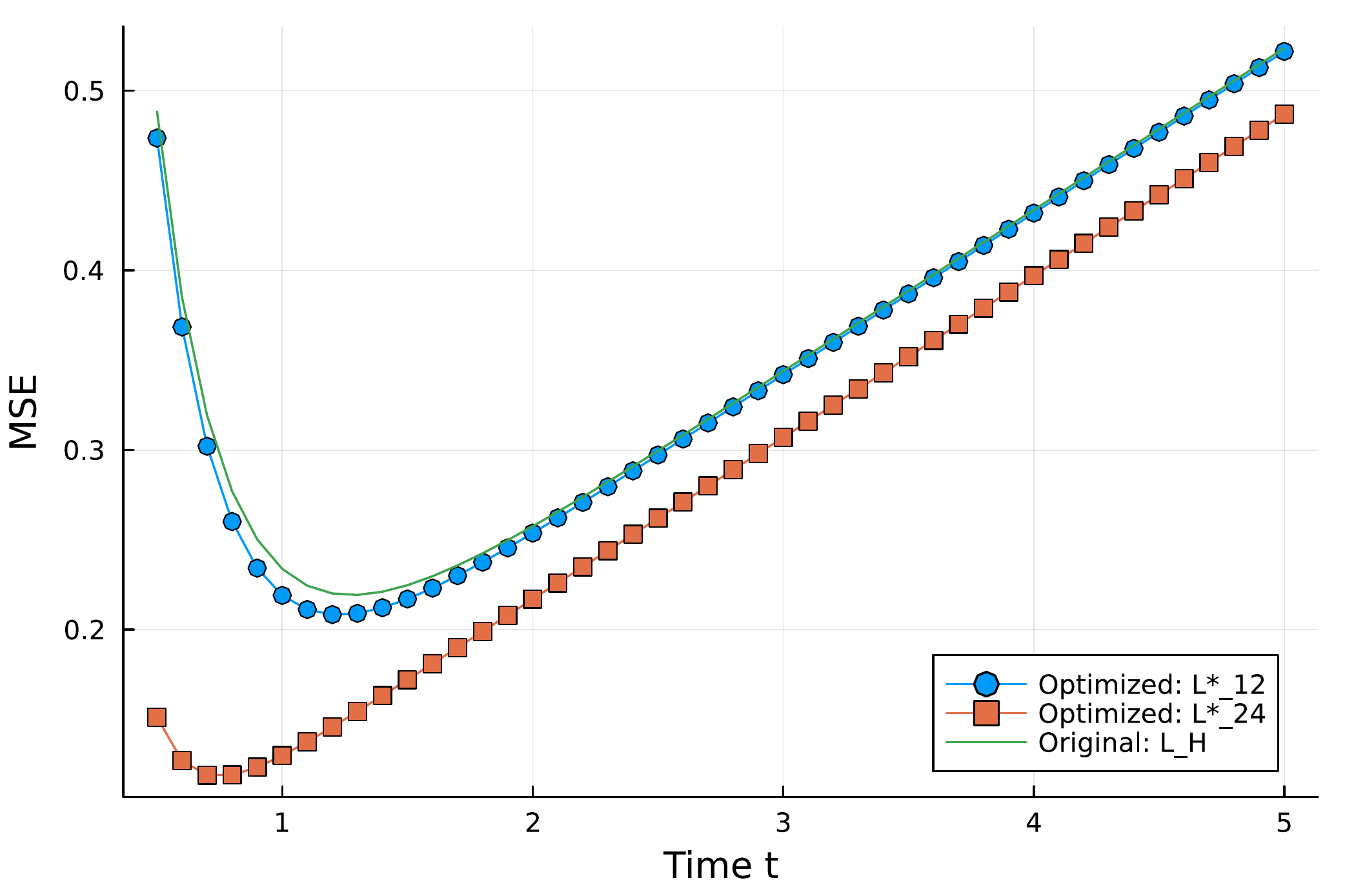}
  \caption{House graph: MSE values of $\bm L^*_{12}, \bm L^*_{20}$ and the unweighted Laplacian matrix $\bm L_H$.}
  \label{MSE_house}
\end{figure}

The MSE values of the optimization result $\bm L^*_{12}, \bm L^*_{24}$ 
and the unweighted Laplacian matrix $\bm L_H$ are shown in Fig.\ref{MSE_house}. 
We can observe that $\bm L_{12}^*$ achieves slightly smaller MSE values 
compared with the unweighted Laplacian matrix $\bm L_H$.
The Laplacian matrix $\bm L^*_{24}$ provides much smaller MSE values
than those of $\bm L_H$.
The sums of inverse eigenvalues are $1.64, 1.59,0.82$ for $\bm L_H$, $\bm L_{12}^*$, and 
$\bm L_{24}^*$, respectively.

\subsection{Barab\'asi-Albert (BA) random graphs (Optimization problem B)}

As an example of random scale-free networks, 
we here handle Barab\'asi-Albert random graph which use a preferential attachment mechanism. 
The number of edges between a new node to existing nodes is assumed to be 5.

In this experiment, we generated an instance of 
Barab\'asi-Albert random graph with 50 nodes. 
The unweighted 
Laplacian of the instance is denoted by $\bm L_B$.
The sum of the diagonal elements of $\bm L_B$ is
$450$. The parameter setting for optimization
is the same as the one used in the previous 
subsection except for $D = 450$.
The output of the optimization algorithm is
referred to as $\bm L^*$.

\begin{figure}[htbp]
  \centering
  \includegraphics[width=0.95\hsize]{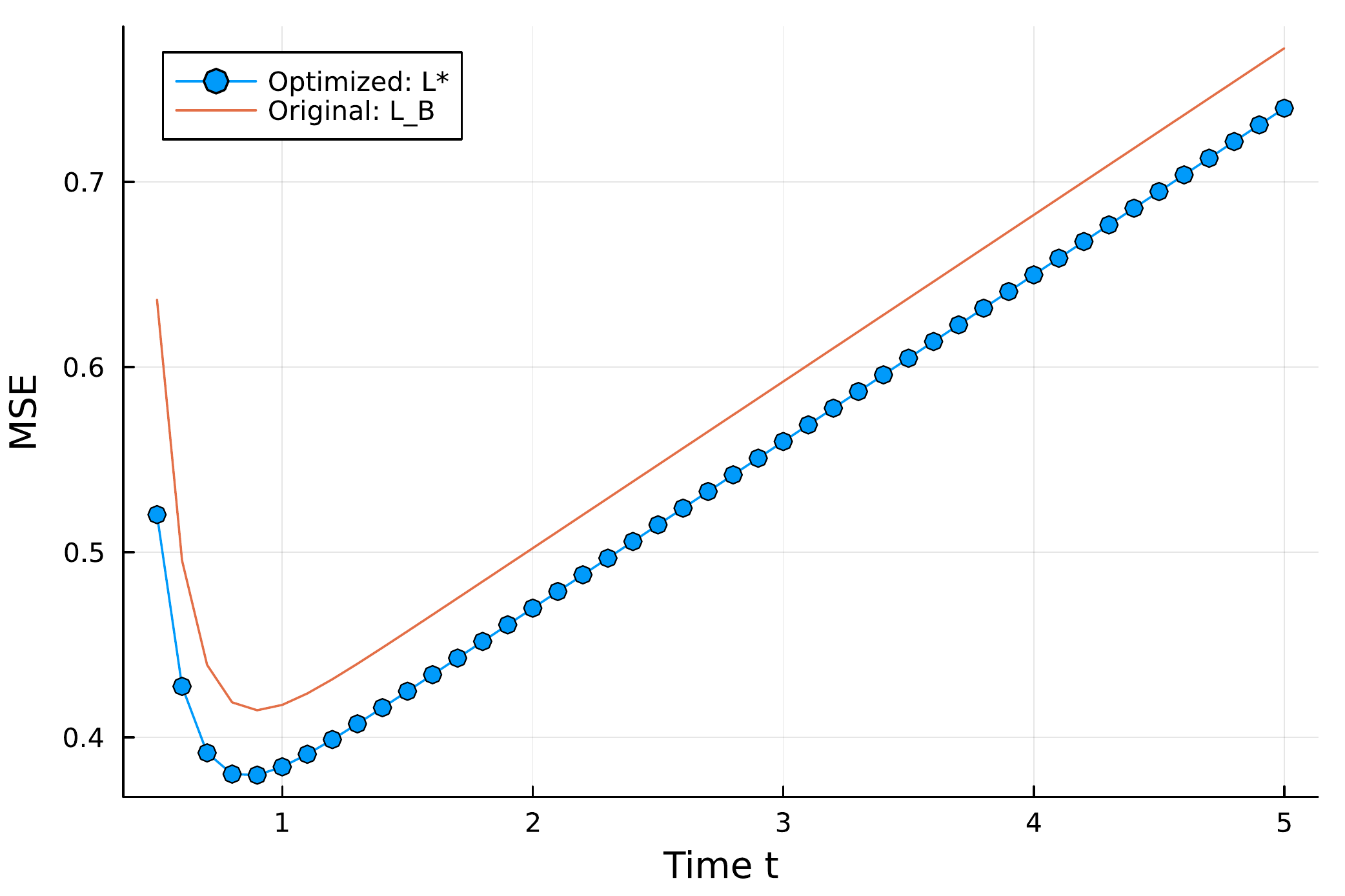}
  \caption{Barab\'asi-Albert random graph: MSE values of $\bm L^*$ and the unweighted Laplacian matrix $\bm L_B$.}
  \label{MSE_BA}
\end{figure}

Figure \ref{MSE_BA} presents the MSE values 
of the original unweighted Laplacian $\bm L_B$ and
the optimization output $\bm L^*$.
We can observe that the optimized MSE values
are substantially smaller than those of the unweighted 
Laplacian $\bm L_B$.
The sums of inverse eigenvalues for $\bm L^*$ and $\bm L_B$ are $6.44$ and $7.16$, respectively.

\begin{figure}[htbp]
  \centering
  \includegraphics[width=0.95\hsize]{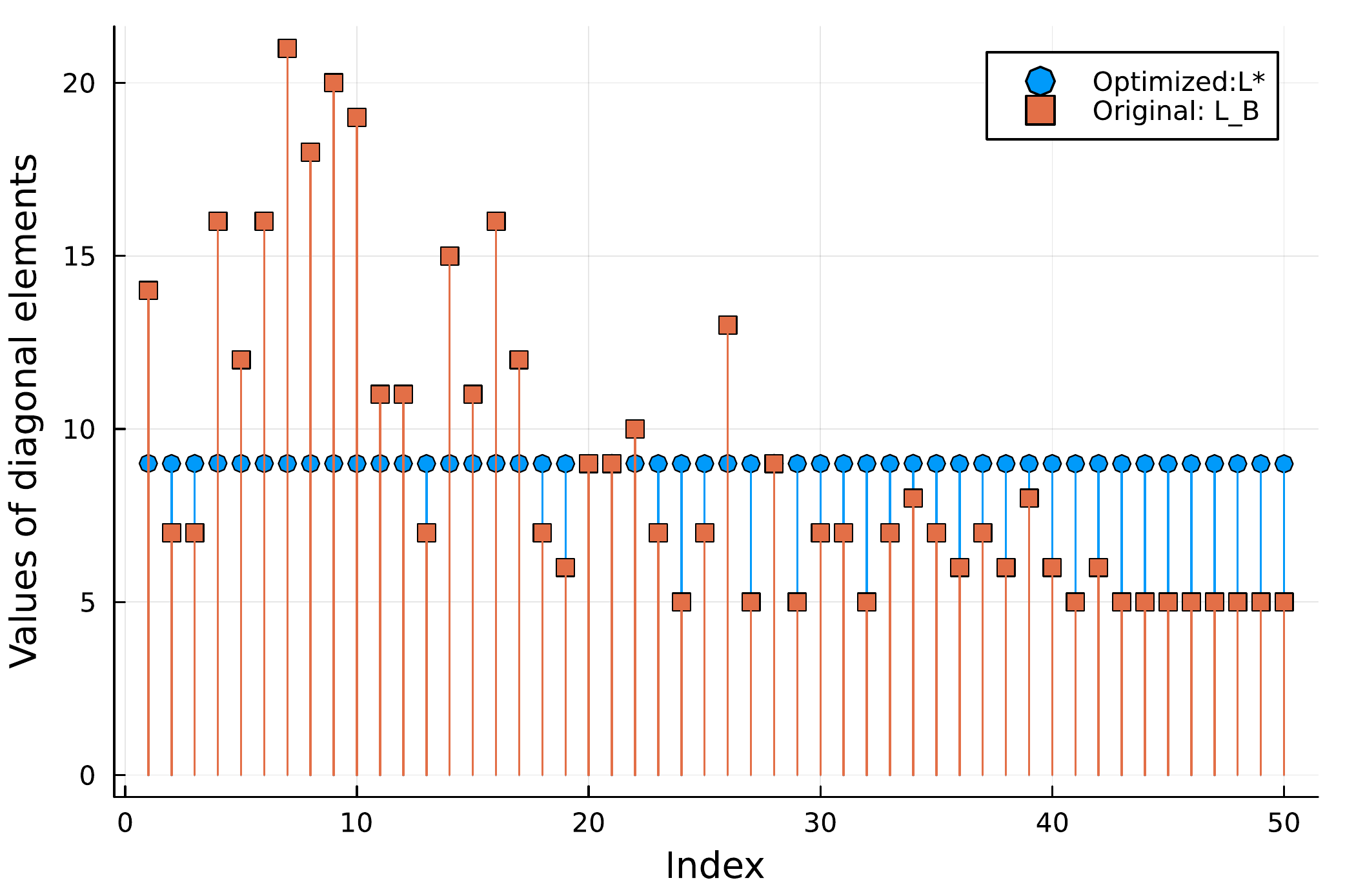}
  \caption{Comparison of diagonal elements of $\bm L^*$ and  $\bm L_B$.}
  \label{degree_BA}
\end{figure}

Figure \ref{degree_BA} illustrates the values of diagonal elements of $\bm L^*$ and  $\bm L_B$.
It can be observed that the values distribution of $\bm L^*$ is almost flat although 
the values of $\bm L_B$ varies from 5 to 21. 
This observation is consistent with the tendency observed in the previous subsection regarding 
the house graph.

\section{Conclusion}
\label{sec:conclusion}

In this paper, we have formulated a noisy average consensus system through a SDE.
This formulation allows for an analytical study of the stochastic dynamics of the system.
We derived a formula for the evolution of covariance for the EM method. 
Through the weak convergence property, we have established Theorem \ref{cov_theorem} 
and derived a MSE formula that provides the MSE at time $t$.
Analysis of the MSE formula reveals that the sum of inverse eigenvalues for the Laplacian matrix 
is the most significant factor impacting the MSE dynamics. 
To optimize the edge weights, a deep unfolding-based technique is presented.
The quality of the solution has been validated by numerical experiments.

It is important to note that the theoretical understanding gained in this study 
will also provide valuable perspective on consensus-based distributed algorithms  
in noisy environments.
In addition, the methodology for optimization proposed in this paper is versatile 
and can be adapted for various algorithms operating on graphs. 
The exploration of potential applications will be an open area for further studies.

 \section*{Acknowledgement}
This study was supported by JSPS Grant-in-Aid for Scientific Research (A) Grant Number 22H00514.
The authors thank Prof. Masaki Ogura for letting us know the related work \cite{Jadbabaie} on 
discrete-time average consensus systems.

\end{document}